\documentclass[12pt]{amsart}
\usepackage{amsmath,amsthm,amsfonts,amssymb,texdraw, eucal}
\numberwithin{equation}{section}
\newtheorem{thm}[equation]{Theorem}

\newtheorem{prop}[equation]{Proposition}
\newtheorem{rem}[equation]{Remark}
\newtheorem{ex}[equation]{Example}
\newtheorem{exs}[equation]{Examples}
\theoremstyle{definition} \newtheorem{defn}[equation]{Definition}
\newcommand{\m}{\medbreak}
\newcommand{\bi}{\bigbreak}

\newcommand{\dotarrow} {\buildrel {* \dots *} \over \longrightarrow}
\newcommand{\dotlarrow} {\buildrel {* \dots *} \over \longleftarrow}
\newcommand{\g}{\mathfrak g}
\newcommand{\h}{\mathfrak h}

\newcommand{\hf}{\frac{1}{2}}

 \newcommand{\iso} {\buildrel \sim \over \rightarrow}
  \newcommand{\iarrow} {\buildrel i \over \rightarrow}
   \newcommand{\iunarrow} {\buildrel {\un {i}} \over \rightarrow}
 \newcommand{\la}{\langle}
 
\newcommand{\ot}{\otimes}
\newcommand{\ov}{\overline}
\newcommand{\p}{\mathfrak p}
\newcommand{\ra}{\rangle}

\newcommand{\un}{\underline}

\newcommand{\zarrow} {\buildrel 0 \over \rightarrow}
\newcommand{\zlarrow} {\buildrel 0 \over \leftarrow}

\begin{document}

\title[PERFECT CRYSTALS]
{\large LEVEL 1 PERFECT CRYSTALS \\ AND PATH REALIZATIONS OF \\
BASIC REPRESENTATIONS AT $\boldsymbol{q=0}$.}

\author{Georgia Benkart$^{\star}$} \address{Department of Mathematics \\ University
of Wisconsin \\  Madison, WI  53706, USA}
\email{benkart@math.wisc.edu}

\author{Igor Frenkel$^{\ast}$}
\address{Department of Mathematics\\
Yale University \\ New Haven, CT 06520, USA }
\email{frenkel@math.yale.edu}

\author{Seok-Jin Kang$^{\dag}$}
\address{Department of Mathematical Sciences and Research Institute of Mathematics\\
Seoul National University \\
San 56-1 Shinrim-dong, Kwanak-ku \\ Seoul 151-747, Korea }
\email{sjkang@math.snu.ac.kr}

\author{Hyeonmi Lee$^{\ddag}$}
\address{School of Mathematics\\
Korea Institute for Advanced Study\\
207-43 Cheongryangri-Dong, Dongdaemun-Gu\\
 Seoul 130-722, Korea  }
\email{hmlee@kias.re.kr}

\thanks{The authors gratefully acknowledge support from the following
grants:\quad $^{\star}$NSF \#{}DMS--0245082 and NSA
MDA904-03-1-0068;    $^{\ast}$NSF \#{}DMS--0070551;  $^{\dag}$KOSEF  \#R01-2003-000-10012-0
and $^{\dag}$KRF \#2003-070-C00001;  and $^{\ddag}$KOSEF \#R01-2003-000-10012-0.}

\date{July 1, 2005} 

\begin{abstract}  We present a uniform construction of level 1 perfect crystals $\mathcal B$
for all affine Lie algebras.   We also introduce the notion of
a crystal algebra and give an explicit description of its multiplication.  
This allows us to determine the energy function on
$\mathcal B \otimes \mathcal B$ completely and thereby give 
a path realization of the basic representations at $q=0$ in the homogeneous
picture.     \end{abstract}
\maketitle

\begin{section}{Introduction} \end{section}  
In the last two decades,  intensive study of solvable lattice models in statistical mechanics
has led to new constructions of representations of affine Lie algebras and their
quantum counterparts.  One noteworthy example is a conjectural realization
of the basic representations for the quantum affine algebra $U_q(\widehat \g)$ for $\g =  {\mathfrak{sl}}_2$  via the states of the XXZ model,

\begin{equation}\label{eq:XXZ}   \quad \cdots \ot V \ot V \ot V,   \end{equation}

\noindent where $V$ is the natural two-dimensional representation of the quantum group
$U_q(\mathfrak {sl}_2)$ extended to a representation of  the quantum affine algebra $U_q'(\widehat {\mathfrak{sl}}_2)$ (the subalgebra of $U_q(\widehat {\mathfrak{sl}}_2)$ without the ``degree operator'').    While its rigorous mathematical meaning is not fully understood for arbitrary values of the parameter $q$, there exists a well-established theory of perfect crystals developed to handle the
limiting case $q = 0$ of this construction.    In particular, the two-dimensional representation
$V$ of $U_q'(\widehat {\mathfrak{sl}}_2)$  in \eqref{eq:XXZ}   admits the structure of the simplest nontrivial perfect crystal of level 1.    This perfect crystal gives rise to
a  path realization of the basic representations of $U_q'(\widehat {\mathfrak{sl}}_2)$ 
and hence to a construction of their crystal graphs at $q = 0$.    Moreover, this construction
has undergone  various generalizations:  \  first to the case  where $V$ is the natural
$n$-dimensional representation of  $U_q'(\widehat{\mathfrak{sl}}_n)$, then to other classical types and to higher level perfect crystals for them
(see \cite{(KMN)$^2$a},  \cite{HK}  and the references therein),    
and lastly  to perfect crystals and path realizations  for some of the exceptional types  (G$_2^{(1)}$ in \cite{Y},
E$_6^{(1)}$ and E$_7^{(1)}$ in \cite{M},  and D$_4^{(3)}$  as announced in \cite{MOY}). 
\medskip

In this paper, we give a uniform construction of level 1 perfect crystals that yields
a path realization of the crystal graphs for the basic representations of  \textbf{all} affine Lie algebras.
Our perfect crystals arise in  the following way.  
Associated to the quantum affine algebra $U_q(\widehat \g)$ is a certain finite-dimensional
simple Lie algebra  $\g$, which is of type  X$_n$ in the untwisted type X$_n^{(1)}$ case  and is given by \eqref{eq:fin}
in the twisted cases.   In our construction,   $V$ as a module for $U_q(\g)$ is a direct
sum  

\begin{equation}\label{eq:modsum}   V = L_{\g}(0) \oplus L_{\g}(\theta)  \end{equation}

\noindent of the trivial one-dimensional $U_q( \g)$-representation
$L_{\g}(0)$ with the representation $L_{\g}(\theta)$  of  $U_q(\g)$
having highest weight $\theta$, which is  the highest root  of 
$ \g$ in the untwisted case and the highest short root of $\g$ in all other
cases but A$_{2n}^{(2)}$.   
Our first result, Theorem 3.5, asserts that the action of $U_q( \g)$ on the module $V$ in
\eqref{eq:modsum}  can be extended to an action of the quantum affine algebra
$U_q'(\widehat \g)$,   and the union

\begin{equation} \label{eq:Bunion}  \mathcal B = \mathcal B(0) \sqcup \mathcal B(\theta)
\end{equation}
of the corresponding crystals has the structure of a perfect crystal of level 1 for $V$.   \medskip

In order to introduce the homogeneous grading in the path realization of the crystal graph of the
basic representations  resulting from our perfect crystal of level 1,  we study in detail the
energy function,

\begin{equation}\label{eq:Hsymbol} H:  \mathcal B \otimes \mathcal B \rightarrow \mathbb Z. \end{equation}
The definition of $H$ implies that its values are constant on each connected component of
the crystal graph $\mathcal B \otimes \mathcal B$ with the 0-arrows omitted.    Thus, the
computation of the energy function \eqref{eq:Hsymbol} easily reduces to specifying its
values on the connected components  of $\mathcal B(\theta) \ot \mathcal B(\theta)$.  We
show that $H = 1$ on all the components except those isomorphic to $\mathcal B(0)$,
$\mathcal B(\theta)$, and $\mathcal B(2\theta)$, where $H = 0,0$, and $2$,  respectively.
We describe explicitly the component of $\mathcal B(\theta) \ot \mathcal B(\theta)$  isomorphic
to $\mathcal B(\theta)$ in Proposition \ref{thetacomp}
and the component isomorphic to $\mathcal B(2\theta)$ in Proposition \ref{2thetacomp}, and thereby 
conclude the determination of the energy function in Theorem \ref{thm:energy}.

\medskip
The description of the component  isomorphic to
$\mathcal B(\theta)$ (there is only one except for type A$_n$,  $n\geq 2$) gives rise to the projection

\begin{equation}\label{eq:multproj}  \mathfrak{m}:  \mathcal B(\theta) \otimes \mathcal B(\theta) \rightarrow \mathcal B(\theta), \end{equation}

\noindent which can be viewed as a multiplication endowing $\mathcal B(\theta)$ with a crystal algebra
structure.   When $\theta$ is the highest root of $ \g$,   then $\mathcal B(\theta)$ is the crystal
graph of the adjoint representation of $ \g$, and \eqref{eq:multproj} yields the structure
of a {\it crystal Lie algebra}.   When $\theta$ is the highest short root for $ \g$, then
the multiplication in \eqref{eq:multproj} gives a crystal version  of the exceptional Jordan algebra for type F$_4$ and
of the octonions for type G$_2$,  without their unit element.   

\medskip
Our construction suggests various generalizations.  First, one can replace the module $V$ in
\eqref{eq:modsum}  by more general finite-dimensional indecomposable modules for
$U_q'(\widehat \g)$
and study when they admit the structure of a perfect crystal.  
In particular, a family of perfect crystals based on Kirillov-Reshetikhin modules was conjectured
in \cite{HKOTY}\footnote{We thank V. Chari for this remark.}.
When a perfect crystal
structure exists, the path realization should provide a relation between the  finite-dimensional and infinite-dimensional representations of quantum affine algebras.   Second, one can try to
``melt the crystals'' and find the quantum versions of the classical identities that
characterize crystal Lie algebras, the crystal exceptional Jordan algebra, and the
crystal octonion algebra.   Finally, one can attempt to make
sense of the tensor product construction of the basic  representations of quantum affine algebras.      

\medskip  Our paper is organized as follows:     Section 2 is devoted to a review of
the notion of a perfect crystal and its corresponding energy function.  In Section 3,
we present our construction of a perfect crystal of level 1 and verify that all the axioms of a perfect crystal are satisfied
except for the connectedness of $\mathcal B \ot \mathcal B$.   The next section  establishes
that $\mathcal B \ot \mathcal B$ is indeed connected.   In Section 5,  we discuss the crystal
algebras that arise from forgetting the 0-arrows and looking
at the connected component of $\mathcal B \ot \mathcal B$
which is isomorphic to $\mathcal B(\theta)$.    In the final section, we evaluate
the energy function on the various components of the crystal $\mathcal B \ot \mathcal B$
minus its  0-arrows.  
\bigskip

\begin{section}{Basics on Perfect Crystals}  \end{section}

In this section, we describe the theory of perfect crystals for affine Lie algebras.
Our discussion here sets the stage for the next part where we give a uniform construction
of  level 1 perfect crystals for all affine Lie algebras.
\m
Let  $\mathcal I = \{0,1,\dots,n\}$ be an index set,  and let $A = \bigl( a_{i,j}\bigr)_{i,j\in \mathcal I}$
be a   Cartan matrix of affine type.  Thus,  $A$ can be characterized by the
following properties:   $a_{i,i} = 2$ for
all $i \in \mathcal I$, \ $a_{i,j} \in \mathbb Z_{\leq 0}$, and $a_{i,j} = 0$ if and only if $a_{j,i} = 0$
for all $i \neq j$ in $\mathcal I$.   The rank of $A$ is $n$,  and if $v  \in \mathbb R^{n+1}$
and $A v \geq 0$ (componentwise), then $v > 0$ or $v = 0$.   We assume $A$ is
indecomposable so that  if $\mathcal I = \mathcal I' \cup \mathcal I''$ where $\mathcal I'$
and $\mathcal I''$ are nonempty, then for some $i \in \mathcal I'$ and $j \in \mathcal I''$,
the entry $a_{i,j} \neq 0$.
An affine Cartan matrix is always symmetrizable -- there exists a diagonal matrix
$D = \hbox{\rm diag}(s_i \mid i \in \mathcal I)$ of positive integers such that
$DA$ is symmetric.  \m

The free abelian group

\begin{equation} Q^\vee = \mathbb Z h_0 \oplus \mathbb Z h_1 \oplus \cdots \oplus \mathbb Z h_n
\oplus \mathbb Z d \end{equation}
is the {\it extended coroot lattice}.  The linear functionals $\alpha_i$ and $\Lambda_i$ ($i \in \mathcal I$)
on the complexification $\h = \mathbb C \ot_{\mathbb Z} Q^\vee$ of $Q^\vee$ given by

\begin{equation}\begin{array}{cccc} \langle h_j,\alpha_i \rangle: =  \alpha_i(h_j) = a_{j,i}& \qquad &\langle d, \alpha_i \rangle: = \alpha_i(d) = \delta_{i,0}& \\
\langle h_j, \Lambda_i \rangle: = \Lambda_i(h_j) = \delta_{i,j}& \qquad &\langle d, \Lambda_i \rangle: = \Lambda_i(d) = 0 &\quad (i,j\in \mathcal I) \end{array}
\end{equation}
are the {\it simple roots} and {\it fundamental weights},  respectively.   Let $\Pi = \{\alpha_i \mid i \in \mathcal I\}$
denote the set of simple roots
and $\Pi^\vee = \{h_i \mid i \in \mathcal I\}$ the set of simple coroots.  The weight lattice
\begin{equation} P = \{\lambda \in \mathfrak h^* \mid \lambda(Q^\vee) \subset \mathbb Z\}\end{equation} contains the set
$P^+ = \{\lambda \in P \mid \lambda(h_i) \in \mathbb Z_{ \geq 0}$ for all $i \in \mathcal I\}$
of {\it dominant integral weights}.  \m

The affine Lie algebra $\widehat \g$ attached to the data $(A,\Pi,\Pi^\vee,P,Q^\vee)$ has generators $e_i,f_i \ (i \in \mathcal I),  h \in \h$,  which satisfy certain relations
(see for example \cite{K} or \cite[Prop.~2.1.6]{HK}).  The algebra $\widehat \g$ can be of three types:

$$\widehat \g= \begin{cases}\hbox{\rm A}_n^{(1)}, \hbox{\rm B}_n^{(1)}, \hbox{\rm C}_n^{(1)}, \hbox{\rm D}_n^{(1)},
\hbox{\rm E}_6^{(1)}, \hbox{\rm E}_7^{(1)}, \hbox{\rm E}_8^{(1)}, \hbox{\rm F}_4^{(1)}, \hbox{\rm G}_2^{(1)}, \\
\hbox{\rm A}_{2n}^{(2)}, \hbox{\rm A}_{2n-1}^{(2)}, \hbox{\rm D}_{n+1}^{(2)}, \hbox{\rm E}_6^{(2)},
  \\
\hbox{\rm D}_4^{(3)}. \end{cases}$$
Associated to   $\widehat \g$ is a finite-dimensional simple Lie algebra $\g$ over $\mathbb C$
given by
\begin{equation}\label{eq:fin} \vbox{
\offinterlineskip \halign{ \strut \vrule  $\hfil#\hfil$  & \vrule
$\hfil#\hfil$  & \vrule $\hfil#\hfil$ & \vrule $\hfil#\hfil$ & \vrule $\hfil#\hfil$ & \vrule  $\hfil#\hfil$  & \vrule  $\hfil#\hfil$
\vrule\cr \noalign{\hrule} \,\widehat \g & \, \hbox{\rm X}_n^{(1)} & \,\hbox{\rm A}_{2n}^{(2)} & \,\hbox{\rm A}_{2n-1}^{(2)}
 & \,\hbox{\rm D}_{n+1}^{(2)} & \,\hbox{\rm E}_6^{(2)} &\, \hbox{\rm D}_4^{(3)} \cr \noalign{\hrule}\,
\g &\, \hbox{\rm X}_n & \,\hbox{\rm C}_n & \, \hbox{\rm C}_n & \,  \hbox{\rm B}_n & \, \hbox{\rm F}_4^t
& \, \hbox{\rm G}_2  \cr \noalign{\hrule}} }\end{equation}
The superscript $t$ is to indicate that the Cartan matrix  is the transpose of the one
used for the $\hbox{\rm F}_4$ associated to $\hbox{\rm F}_4^{(1)}$.    \m

The {\it canonical central element} $c$ and the {\it null root} $\delta$  are given by the
expressions

\begin{equation}\label{eq:centnull}\begin{array}{cc} &c = c_0 h_0 + c_1 h_1 + \cdots + c_n h_n, \\
&\delta = d_0 \alpha_0 + d_1 \alpha_1 + \cdots + d_n \alpha_n,
\end{array}\end{equation}

\noindent where $c_0 = 1$,   and
$d_0 = 1$  except for type A$_{2n}^{(2)}$ where $d_0 = 2$.
The first term comes from the fact that the center of the corresponding affine Lie algebra $\widehat \g$ is  generated by $c$,   while the second comes from the fact that the vector $[d_0,d_1,\dots, d_n]^{t}
\in \mathbb C^{n+1}$ spans the null space of the Cartan matrix $A$.
We say that a dominant weight $\lambda \in P^+$ has {\it level \  $\ell$} \ if
$\langle c, \lambda \rangle : = \lambda(c) = \ell$.  \m

Given any $n \in \mathbb Z$ and an indeterminate $x$,  let
$$[n]_x  = \frac{x^n- x^{-n}}{x-x^{-1}}.$$
Set $[0]_x\,! = 1$ and $[n]_x\,! = [n]_x[n-1]_x \cdots [1]_x$ for $n \geq 1$,  and for $m \geq n \geq 0$,   let
$$\left[\begin{array}{c} m \\ n \end{array}\right]_{x} = \frac{[m]_x\,!}{[n]_x\,!\,[m-n]_x\,!}.$$

  \bigbreak

\begin{defn}\label{defn:qalg}The {\it quantum affine algebra} $U_q(\widehat \g)$ associated with $(A,\Pi,\Pi^\vee, P, Q^\vee)$
is the associative algebra with unit element over $\mathbb C(q)$ (where $q$ is an indeterminate)
with generators $e_i,\,f_i \ (i \in \mathcal I)$,  $q^h \ (h \in Q^\vee)$  satisfying the defining relations
\begin{itemize}
\item[\rm{(1)}] $ q^0 = 1, \ q^hq^{h'} = q^{h+h'}$  \ \ \  for $h,h' \in Q^\vee$,
\item[\rm{(2)}] $q^h e_i q^{-h} = q^{\alpha_i(h)}e_i$ \ \ \ for $h \in Q^\vee$, $i \in \mathcal I$,
\item[\rm{(3)}] $q^h f_i q^{-h} = q^{-\alpha_i(h)}f_i$ \ \ \  for $h \in Q^\vee$, $i \in \mathcal I$,
\item[\rm{(4)}] $e_if_j - f_je_i = \displaystyle{\delta_{i,j}\frac{K_i-K_i^{-1}}{q_i-q_i^{-1}}}$ \ \ \  for $i,j \in \mathcal I$,
\item[\rm{(5)}] $\displaystyle{\sum_{k=0}^{1-a_{i,j}} \left[\begin{array}{c} 1-a_{i,j} \\ k \end{array}\right]_{q_i}
e_i^{1-a_{i,j}-k} e_j e_i^k = 0}$  \ \ \  for  $i \neq j$,
\item[\rm{(6)}] $\displaystyle{\sum_{k=0}^{1-a_{i,j}} \left[\begin{array}{c} 1-a_{i,j} \\ k \end{array}\right]_{q_i}
f_i^{1-a_{i,j}-k} f_j f_i^k = 0}$  \ \ \ for  $i \neq j$,
 \end{itemize}
where $q_i = q^{s_i}$ and $K_i = q^{s_i h_i}$.
\end{defn}

Crystal base theory  has been developed for $U_q(\widehat \g)$-modules in the category
$\mathcal O_{\hbox {\rm \small  int}}$ of integrable modules.    This is the category of $U_q(\widehat \g)$-modules
$M$ such that \smallskip
\begin{itemize}
\item[{(a)}] $M$ has a weight space decomposition: \
$M = \bigoplus_{\lambda \in P} M_\lambda$,  where
$M_\lambda = \{v \in M \mid q^h.v = q^{\lambda(h)}v$ for all $h \in Q^\vee \}$;
\item[{(b)}] there are finitely many $\lambda_1,\dots, \lambda_k \in P$ such that
wt$(M) \subseteq \Omega(\lambda_1) \cup \dots \cup \Omega(\lambda_k)$,\
where wt$(M) = \{\lambda \in P \mid M_\lambda \neq 0\}$ and
$\Omega(\lambda_j) = \{\mu \in P \mid \mu \in \lambda_j+ \sum_{i\in \mathcal I} \mathbb Z_{\leq 0} \alpha_i\}$;
\item[{(c)}]  the elements $e_i$ and $f_i$ act locally nilpotently on $M$ for all $i \in \mathcal I$.
\end{itemize} \m

If $M$ is a module in category $\mathcal O_{\hbox {\rm \small  int}}$, then for each $i \in \mathcal I$,  a weight vector $u \in M_\lambda$
has a unique expression $u = \sum_{k=0}^N  f_i^{(k)} u_k$,  where
$u_k \in M_{\lambda+k\alpha_i} \cap \ker e_i$ for   $k=0,1,\dots,N$,  and
$f_i^{(k)} = f_i^k/ [k]_{q_i}!$.    The {\it Kashiwara operators}
are defined on $u$ using these expressions according to the rules

\begin{equation} \tilde e_i u = \sum_{k=1}^N f_i^{(k-1)}u_k,   \qquad \ \   \tilde f_i u = \sum_{k=0}^N
f_i^{(k+1)}u_k.  \end{equation}

Let $\mathbb A_0 = \{ f = g/h \mid g,h  \in \mathbb C[q], \  h(0) \neq 0 \}$ be the localization of $\mathbb C[q]$
at the ideal $(q)$.   Every module $M \in \mathcal O_{\hbox {\rm \small  int}}$ has a special type
of $\mathbb A_0$-lattice called a crystal lattice.

\begin{defn}  Assume $M$ is a $U_q(\widehat \g)$-module in category $\mathcal O_{\hbox {\rm \small  int}}$.  A free
$\mathbb A_0$-submodule $\mathcal L$ of $M$ is a {\it crystal lattice} if
\begin{itemize}
\item[{\rm (i)}] $\mathcal L$ generates $M$ as a vector space over $\mathbb C(q)$;
\item[{\rm (ii)}]  $\mathcal L = \bigoplus_{\lambda \in P} \mathcal L_\lambda$ \  where \
$\mathcal L_\lambda = M_\lambda \cap \mathcal L$;
\item[{\rm (iii)}]  $\tilde e_i \mathcal L \subset \mathcal L$ \  and \
$\tilde f_i \mathcal L \subset \mathcal L$.
\end{itemize}  \end{defn}

Since the operators $\tilde e_i$ and $\tilde f_i$ preserve the lattice $\mathcal L$,
they give well-defined operators on the quotient $\mathcal L/q \mathcal L$, which we denote by
the same symbols.

\begin{defn}  A {\it crystal base} for a $U_q(\widehat \g)$-module $M \in \mathcal O_{\hbox {\rm \small  int}}$ is a pair
$(\mathcal L, \mathcal B)$   such that
\begin{itemize}
\item[{\rm (1)}] $\mathcal L$ is a crystal lattice of $M$;
\item[{\rm (2)}]  $\mathcal B$ is a $\mathbb C$-basis of $\mathcal L/q \mathcal L \cong \mathbb C \ot_{\mathbb A_0} \mathcal L$;
\item[{\rm(3)}]  $\mathcal B = \sqcup_{\lambda \in P} \mathcal B_\lambda$, \ where \
$\mathcal B_\lambda = \mathcal B \cap (\mathcal L_\lambda/ q \mathcal L_\lambda)$;
\item[{\rm(4)}]  $\tilde e_i \mathcal B \subset \mathcal B \cup \{0\}$ \  and \
$\tilde f_i \mathcal B \subset \mathcal B \cup \{0\}$ for all $i \in \mathcal I$;
\item[{\rm(5)}]  $\tilde f_i b = b'$ \  if and only if \  $b = \tilde e_i b'$ \  for \  $b,b' \in \mathcal B$ \ and \ $i\in \mathcal I$.  \end{itemize}
\end{defn}

As each $M \in \mathcal O_{\hbox {\rm \small  int}}$ has such a crystal base, we can associate to each $M$ a
crystal graph having $\mathcal B$  as the set of vertices.  Vertices $b,b' \in \mathcal B$
are connected by an arrow labelled by $i$ pointing from $b$ to $b'$  if and only if $\tilde f_i b =
b'$.   The crystal graph encodes much of the combinatorial information about $M$.
\m

For $i \in \mathcal I$,  let \  $\varepsilon_i, \varphi_i: \mathcal B \rightarrow \mathbb Z$ \ be defined
by
\begin{equation}
\begin{array}{cc}
&\varepsilon_i(b) = \max\{ k  \geq 0 \mid \tilde e_i^k b \in \mathcal B\}, \\
&\varphi_i(b) =  \max\{ k  \geq 0 \mid \tilde f_i^k b \in \mathcal B\}. \end{array}
\end{equation}
{F}rom property (5) we see that $\varepsilon_i(b)$ is just the number of $i$-arrows
coming into $b$ in the crystal graph and $\varphi_i(b)$ is just the number of
$i$-arrows emanating from $b$.    Moreover, $\varphi_i(b) -\varepsilon_i(b) = \lambda(h_i)$ for
all $b \in \mathcal B_\lambda$.    Thus, if

$$\varepsilon(b) = \sum_{i\in \mathcal I} \varepsilon_i(b) \Lambda_i,
\qquad \varphi(b) = \sum_{i \in \mathcal I} \varphi_i(b) \Lambda_i, $$
then \ $\hbox{\rm wt}\,b = \varphi(b) -\varepsilon(b) = \lambda$ \  for all $b \in \mathcal B_\lambda$. \m

Morphisms between two crystals \, $\mathcal B_1$ \, and \,
$\mathcal B_2$ \, associated to \, $U_q(\widehat \g)$ \,  are maps
$\Psi: \mathcal B_1 \cup \{0\} \rightarrow \mathcal B_2 \cup
\{0\}$ such that  $\Psi(0) = 0$; \  $\hbox{\rm wt}\,\Psi(b) =
\hbox{\rm wt}\,b$, \ $\varepsilon_i(\Psi(b)) = \varepsilon_i(b)$
and  $\varphi_i(\Psi(b)) = \varphi_i(b)$ for all $i \in \mathcal I$;  and if $b, b' \in \mathcal B_1$,
$\tilde f_i b = b'$,  and $\Psi(b), \Psi(b') \in \mathcal B_2$, then $\tilde
f_i \Psi(b)= \Psi(b')$, $\tilde e_i \Psi(b')=\Psi(b)$.    A morphism
$\Psi$ is called {\it strict} if it commutes with the Kashiwara
operators $\tilde e_i$, $\tilde f_i$ for all $i\in \mathcal I$.
 
One of the most striking features of crystal bases is their behavior under tensor products.
If $M_j \in \mathcal O_{\hbox {\rm \small  int}}$ for $j=1,2$ and $(\mathcal L_j,\mathcal B_j)$ are the corresponding
crystal bases,  set $\mathcal L = \mathcal L_1 \ot_{{\mathbb A}_0} \mathcal L_2$ and $\mathcal B
= \mathcal B_1 \times \mathcal B_2$.   Then $(\mathcal L, \mathcal B)$ is a crystal base
of $M_1 \ot_{\mathbb C(q)} M_2$, where the action of the Kashiwara operators on $\mathcal B$
is given by

\begin{equation}\begin{array}{cc}
&\tilde e_i (b_1 \ot b_2) = \begin{cases} \tilde e_i b_1 \ot b_2 \quad & \hbox{\rm if} \ \ \varphi_i(b_1)  \geq \varepsilon_i(b_2), \\
b_1 \ot \tilde e_i b_2 \quad & \hbox{\rm if} \ \ \varphi_i(b_1) < \varepsilon_i(b_2), \end{cases}\\
&\tilde f_i (b_1 \ot b_2) = \begin{cases} \tilde f_i b_1 \ot b_2 \quad & \hbox{\rm if} \ \ \varphi_i(b_1) > \varepsilon_i(b_2), \\
b_1 \ot \tilde f_i b_2 \quad & \hbox{\rm if} \ \ \varphi_i(b_1) \leq \varepsilon_i(b_2). \end{cases}
\end{array} \end{equation}

\m

Corresponding to any $\lambda \in P$ is
a one-dimensional module $\mathbb C v_\lambda$ for
the subalgebra  $U^{ \geq 0}$ of $U_q(\widehat \g)$ generated
by $e_i \ (i \in \mathcal I$),  $q^h \ (h \in Q^\vee)$,   where the $U^{ \geq 0}$-action
is given by $e_i.v_\lambda = 0$ and $q^h.v_\lambda = q^{\lambda(h)}v_\lambda$.
The induced module $V(\lambda):= U_q(\widehat \g) \ot_{U^{ \geq 0}} \mathbb C v_\lambda$
(the so-called Verma module) has a unique maximal submodule and a unique
irreducible quotient $L(\lambda)$.      The modules
$L(\lambda)$ for $\lambda \in P^+$ account for
all the irreducible modules in category $\mathcal O_{\hbox {\rm \small int}}$.
Let $(\mathcal L(\lambda), \mathcal B(\lambda))$ denote the crystal
base corresponding to $L(\lambda)$.  Since the
weight space of $L(\lambda)$ corresponding to the
weight $\lambda$ is one-dimensional, the crystal $\mathcal B(\lambda)$
has a unique element of weight $\lambda$, which we denote by $u_\lambda$
in the sequel.   It has the property that $\tilde e_i u_\lambda = 0$
for all $i \in \mathcal I$.   \m

The subalgebra $U_q'(\widehat \g)$ of $U_q(\widehat \g)$ generated by $e_i, f_i, K_i^{\pm 1}\ (i \in \mathcal I)$
is often also referred to as the quantum affine algebra.  The main difference between
$U_q(\widehat \g)$ and $U_q'(\widehat \g)$ is that $U_q'(\widehat \g)$ admits  nontrivial  finite-dimensional irreducible
modules, while $U_q(\widehat \g)$  does not.    The theory of perfect crystals, which we introduce
next, requires us to work with $U_q'(\widehat \g)$ for that reason.
Here we will need the coroot lattice

\begin{equation}\bar Q^\vee = \mathbb Zh_0 \oplus \mathbb Z h_1 \oplus \cdots \oplus \mathbb Z h_n,
\end{equation}
 and its complexification $\bar \h = \mathbb C \ot_{\mathbb Z}\bar Q^\vee$.  When elements of the
 $\mathbb Z$-submodule

\begin{equation}\bar P = \mathbb Z \Lambda_0 \oplus \mathbb Z \Lambda_1 \oplus \cdots \oplus \mathbb Z \Lambda_n  \end{equation}
of  $P$ are restricted to $\bar Q^\vee$, they give the dual lattice of  {\it classical weights}.     Let
  \,$\bar P^+
: = \sum_{i= 0}^n  \mathbb Z_{ \geq 0} \Lambda_i$ \, denote the corresponding set of  dominant weights.  \m

 Every symmetrizable Kac-Moody Lie algebra has a crystal base theory. In particular,
  the finite-dimensional simple Lie algebras ${\g}$ over $\mathbb C$ have a crystal base
  theory, and every finite-dimensional ${\g}$-module has a crystal base.  We refer
  to such crystal bases as {\it finite classical crystals}.   \m

We recall the definition of a perfect crystal (see for example  \cite[Defn.~10.5.1]{HK}).

\m
\begin{defn}\label{defn:pc}   For a positive integer \,$\ell$, we say a finite
classical crystal \,$\mathcal B$ \,  is a \  {\it perfect crystal of
level \,$\ell$} \  for the quantum affine algebra
\,$U_q(\widehat \g)$\,  if
\begin{itemize}
\item[{\rm (1)}]  there is a finite-dimensional $U_q'(\widehat \g)$-module
with a crystal base whose crystal graph is isomorphic to
\, $\mathcal B$ (when the $0$-arrows are removed);
\item[{\rm (2)}]  $\mathcal B \otimes \mathcal B$ \,  is connected;
\item[{\rm (3)}]  there exists a classical weight \,$\lambda_0$\,
such that

$$\hbox{\rm wt}(\mathcal B) \subset \lambda_0 + \frac{1}{d_0}\sum_{i \neq 0} \mathbb Z_{\leq 0} \alpha_i
\quad \hbox{\rm and} \quad  |\mathcal B_{\lambda_0}| = 1;$$

\item[{\rm (4)}]  for any $b \in \mathcal B$, we have \, $\langle c,\varepsilon(b)\rangle
= \sum_{i \in \mathcal I} \varepsilon_i(b) \Lambda_i(c)  \geq \ell$;

\item[{\rm (5)}]   for each \, $\lambda \in \bar P_\ell^+ := \{ \mu \in \bar P^+  \mid \langle c, \mu \rangle = \ell\}$, \, there exist unique
vectors \,$b^\lambda$\, and \,$b_\lambda$\,  in $\mathcal B$ such that
\,$\varepsilon(b^\lambda) = \lambda$\, and \,$\varphi(b_\lambda) = \lambda$.
\end{itemize}
\end{defn}

The significance of perfect crystals is that they provide a means of constructing
the crystal base $\mathcal B(\lambda)$ of any irreducible $U_q(\widehat \g)$-module
$L(\lambda)$ corresponding to a classical weight  $\lambda \in \bar P^+$.    \bi

\begin{thm}\cite{(KMN)$^2$a}  Assume $\mathcal B$ is a perfect crystal of level $\ell > 0$.
Then for any classical weight  $\lambda \in \bar P^+_\ell$,  there is
a crystal isomorphism

\begin{eqnarray*} \Psi:  \mathcal B(\lambda) &\,\,\iso \ \, \mathcal B(\varepsilon(b_\lambda)) \ot \mathcal B \\
u_\lambda&\mapsto \ \ u_{\varepsilon(b_\lambda)}  \ot b_\lambda.
\end{eqnarray*} \end{thm}

\bi

As a consequence of this theorem, any $\lambda \in \bar P^+_\ell$ gives rise to  a sequence of weights and corresponding  elements in  the perfect crystal $\mathcal B$,

\begin{equation}\label{eq:lamb}\begin{array} {ccc}
\lambda_0 = \lambda &\qquad b_0 = b_\lambda &  \\
\lambda_{k+1} = \varepsilon(b_{\lambda_k})
 &\qquad \, b_{k+1} = b_{\lambda_{k+1}}  & \qquad  \hbox{\rm for all}\ \
k  \geq 1,   \end{array} \end{equation}
such that

\begin{gather} \mathcal B(\lambda_k) \ \iso \  \mathcal B(\lambda_{k+1}) \ot \mathcal B    \\
  \quad u_{\lambda_k} \, \mapsto \ \ u_{\lambda_{k+1}}  \ot b_k.  \nonumber \end{gather}
 Iterating this isomorphism, we have
 \begin{gather}
  \mathcal B(\lambda)\ \iso\ \mathcal B(\lambda_1) \ot \mathcal B \ \iso \  \mathcal B(\lambda_2) \ot \mathcal B \ot \mathcal B\ \iso \cdots    \\  
  u_\lambda\  \ \ \mapsto \ \ \ u_{\lambda_1} \ot b_0 \ \ \ \mapsto \ \ \ u_{\lambda_2} \ot b_1 \ot b_0 \ \  \ \mapsto 
  \cdots  \nonumber  \end{gather}

  \bi
  \begin{defn}  For $\lambda \in \bar P^+_\ell$,   the {\it ground state path of weight} $\lambda$ is the tensor product

 $${\p}_\lambda = \,\bigl(b_k)_{k=0}^\infty \,= \  \  \cdots \ot b_{k+1} \ot b_k \ot \cdots \ot b_1 \ot b_0,$$
 where the elements $b_k \in \mathcal B$ are  as in \eqref{eq:lamb}.
 A tensor product $\p = (p_k)_{k=0}^\infty =  \cdots \ot p_{k+1} \ot p_k \ot \cdots \ot p_1 \ot p_0$ of elements $p_k \in \mathcal B$ is said to be a $\lambda$-{\it path} if
 $p_k = b_k$ for all $k \gg  0$.   \end{defn}
 \bi

 \begin{thm}\label{thm:cryiso}\cite{(KMN)$^2$a}  Assume $\lambda \in  \bar P^+_\ell$.
 Then there is a crystal isomorphism
 $\mathcal B(\lambda)  \iso \mathcal P(\lambda)$, with $u_\lambda \mapsto \p_\lambda$,
 between the crystal base $\mathcal B(\lambda)$ of $L(\lambda)$ and the set $\mathcal P(\lambda)$ of $\lambda$-paths.
 \end{thm}  \m

 The crystal structure of $\mathcal P(\lambda)$ referred to in Theorem \ref{thm:cryiso} may be described
 as follows.   Given any $\p = (p_k)_{k=0}^\infty \in \mathcal P(\lambda)$, let $N > 0$ be such that  $p_k = b_k$ for all $k  \geq N$.    As in \cite[(10.48)]{HK}, set

\begin{eqnarray}\label{eq:path}  
&&\ov {\hbox{\rm wt} \p} = \lambda_N + \sum_{k=0}^{N-1} \ov{\hbox{\rm wt} p}_k, \nonumber \\
&&\tilde e_i \p\  = \ \ \cdots \ot p_{N+1} \ot \tilde e_i\left( p_N \ot \cdots \ot p_0 \right), \nonumber \\
&&\tilde f_i \p \ = \ \ \cdots \ot p_{N+1} \ot \tilde f_i\left( p_N \ot \cdots \ot p_0 \right), \\
&&\varepsilon_i(\p) = \max\left( \varepsilon_i(\p')-\varphi_i(b_N), 0\right), \nonumber  \\
&&\varphi_i(\p) = \varphi_i(\p') +  \max\left(\varphi_i(b_N)-\varepsilon_i(\p'), 0\right), \nonumber
\end{eqnarray} 
where $\p' := p_{N-1} \ot \cdots \ot p_1 \ot p_0$ and
$\ov{\hbox{\rm wt}}$ signifies the classical weight of an element
of $\mathcal B$ or $\mathcal P(\lambda)$. \m
 
Equation \eqref{eq:path} describes  the classical weight
$\ov{\hbox{\rm wt} \p}$ (i.e. the element of $\bar P$  attached to each $\p \in
\mathcal P(\lambda)$).     We would like to calculate the actual affine weight  $\hbox{\rm wt}\p$
in $P$.   For this, we need the notion of an energy function. \m

\begin{defn}  Let $V$ be a finite-dimensional $U_q'(\widehat \g)$-module with crystal base
$(\mathcal L, \mathcal B)$.  An {\it energy function} on $\mathcal B$ is a map $H: \mathcal B \ot \mathcal B \rightarrow
\mathbb Z$ satisfying

\begin{equation}\label{eq:ef} H\left(\tilde e_i (b_1 \ot b_2)\right)
= \begin{cases} H(b_1 \ot b_2) & \qquad \hbox{\rm if} \ \ i \neq 0, \\
H(b_1 \ot b_2) + 1 & \qquad \hbox{\rm if} \ \ i = 0 \ \hbox{\rm and} \ \varphi_0(b_1)  \geq \varepsilon_0(b_2)  \\
H(b_1 \ot b_2) - 1 & \qquad \hbox{\rm if} \ \ i = 0 \ \hbox{\rm and} \ \varphi_0(b_1)< \varepsilon_0(b_2).
\end{cases} \end{equation}
for all $b_1,b_2 \in \mathcal B$ with $\tilde e(b_1 \ot b_2) \in \mathcal B \ot \mathcal B$.
\end{defn}
\bi

\begin{ex}  Let $\widehat \g$ be the affine Lie algebra \ $\hbox{\rm A}_2^{(1)}$,  and let
$\mathcal B$ be the crystal with 3 elements:    \end{ex}

\begin{center}
\begin{texdraw}%
\drawdim in
\arrowheadsize l:0.065 w:0.03
\arrowheadtype t:F
\fontsize{7}{7}\selectfont
\textref h:C v:C
\drawdim em
\setunitscale 1.9
\move(2 0)

\move(-1 -1)
\bsegment
\htext(-1.3 0){\normalsize$\mathcal B$}
\htext(-0.8 0){\textup{:}}
\move(1 -0.4)
\bsegment
\move(0 0)\rlvec(0.8 0)\rlvec(0 0.8)\rlvec(-0.8 0)\rlvec(0 -0.8)
\htext(0.38 0.38){$1$}
\esegment
\move(2.02 0)\ravec(0.8 0)
\move(4.02 0)\ravec(0.8 0)
 \move (5.3 -0.5)
 \clvec (5.3 -1.4)(1.5 -1.4)(1.5 -0.5)
\htext(3.3 -1.5){$0$}
\htext(4.3 0.3){$2$}
\htext(2.3 0.3){$1$}
\move(3 -0.4)
\bsegment
\move(0 0)\rlvec(0.8 0)\rlvec(0 0.8)\rlvec(-0.8 0)\rlvec(0 -0.8)
\htext(0.38 0.38){$2$}
\esegment
\move(5 -0.4)
\bsegment
\move(0 0)\rlvec(0.8 0)\rlvec(0 0.8)\rlvec(-0.8 0)\rlvec(0 -0.8)
\htext(0.38 0.38){$3$}
\esegment
\move(1.5 -0.55)\ravec(-0.1 0.2)
\esegment

\move(0 -1)
\bsegment
\htext(-3 -2.6){\normalsize$\mathcal B\otimes\mathcal B$}
\htext(-1.8 -2.6){\textup{:}}
\move(0 0)
\bsegment
\move(0 -3)
\bsegment
\move(0 0)\rlvec(0.8 0)\rlvec(0 0.8)\rlvec(-0.8 0)\rlvec(0 -0.8)
\htext(0.38 0.38){$1$}
\htext(1.15 0.4){$\otimes$}
\esegment
\move(1.5 -3)
\bsegment
\move(0 0)\rlvec(0.8 0)\rlvec(0 0.8)\rlvec(-0.8 0)\rlvec(0 -0.8)
\htext(0.38 0.38){$1$}
\move(1.4 0.4)\ravec(1.5 0)
 \move (-1.9 0.2)
 \clvec (-2.8 0.2)(-2.8 -5.4)(-1.9 -5.4)
\move(-1.95 0.18)\ravec(0.2 0.2)
\htext(-2.9 -2.65){$0$}
\htext(2 0.7){$1$}
\esegment
\esegment

\move(5 0)
\bsegment
\move(0 -3)
\bsegment
\move(0 0)\rlvec(0.8 0)\rlvec(0 0.8)\rlvec(-0.8 0)\rlvec(0 -0.8)
\htext(0.38 0.38){$2$}
\htext(1.15 0.4){$\otimes$}
\esegment
\move(1.5 -3)
\bsegment
\move(0 0)\rlvec(0.8 0)\rlvec(0 0.8)\rlvec(-0.8 0)\rlvec(0 -0.8)
\htext(0.38 0.38){$1$}
\move(1.4 0.4)\ravec(1.5 0)
\move(-0.35 -0.5)\ravec(0 -1.2)
 \move (-1.9 0.2)
 \clvec (-2.8 0.2)(-2.8 -5.4)(-1.9 -5.4)
\move(-1.95 0.18)\ravec(0.2 0.2)
\htext(-2.9 -2.65){$0$}
\htext(2 0.7){$2$}
\htext(-0.7 -0.9){$1$}
\esegment
\esegment

\move(10 0)
\bsegment
\move(0 -3)
\bsegment
\move(0 0)\rlvec(0.8 0)\rlvec(0 0.8)\rlvec(-0.8 0)\rlvec(0 -0.8)
\htext(0.38 0.38){$3$}
\htext(1.15 0.4){$\otimes$}
\esegment
\move(1.5 -3)
\bsegment
\move(0 0)\rlvec(0.8 0)\rlvec(0 0.8)\rlvec(-0.8 0)\rlvec(0 -0.8)
\htext(0.38 0.38){$1$}
\move(-0.35 -0.5)\ravec(0 -1.2)
\htext(-0.7 -0.9){$1$}
\esegment
\esegment

\move(0 -3)
\bsegment
\move(0 -3)
\bsegment
\move(0 0)\rlvec(0.8 0)\rlvec(0 0.8)\rlvec(-0.8 0)\rlvec(0 -0.8)
\htext(0.38 0.38){$1$}
\htext(1.15 0.4){$\otimes$}
\esegment
\move(1.5 -3)
\bsegment
\move(0 0)\rlvec(0.8 0)\rlvec(0 0.8)\rlvec(-0.8 0)\rlvec(0 -0.8)
\htext(0.38 0.38){$2$}
\move(-0.35 -0.5)\ravec(0 -1.2)
 \move (9.4 -0.5)
 \clvec (9.4 -1.5)(0 -1.5)(0 -0.5)
\move(0 -0.55)\ravec(-0.1 0.2)
\htext(4.56 -4.65){$0$}
\htext(-0.7 -0.9){$2$}
\esegment
\esegment

\move(5 -3)
\bsegment
\move(0 -3)
\bsegment
\move(0 0)\rlvec(0.8 0)\rlvec(0 0.8)\rlvec(-0.8 0)\rlvec(0 -0.8)
\htext(0.38 0.38){$2$}
\htext(1.15 0.4){$\otimes$}
\esegment
\move(1.5 -3)
\bsegment
\move(0 0)\rlvec(0.8 0)\rlvec(0 0.8)\rlvec(-0.8 0)\rlvec(0 -0.8)
\htext(0.38 0.38){$2$}
\move(1.4 0.4)\ravec(1.5 0)
\htext(2 0.7){$2$}
\esegment
\esegment

\move(10 -3)
\bsegment
\move(0 -3)
\bsegment
\move(0 0)\rlvec(0.8 0)\rlvec(0 0.8)\rlvec(-0.8 0)\rlvec(0 -0.8)
\htext(0.38 0.38){$3$}
\htext(1.15 0.4){$\otimes$}
\esegment
\move(1.5 -3)
\bsegment
\move(0 0)\rlvec(0.8 0)\rlvec(0 0.8)\rlvec(-0.8 0)\rlvec(0 -0.8)
\htext(0.38 0.38){$2$}
\move(-0.35 -0.5)\ravec(0 -1.2)
\htext(-0.7 -0.9){$2$}
\esegment
\esegment

\move(0 -6)
\bsegment
\move(0 -3)
\bsegment
\move(0 0)\rlvec(0.8 0)\rlvec(0 0.8)\rlvec(-0.8 0)\rlvec(0 -0.8)
\htext(0.38 0.38){$1$}
\htext(1.15 0.4){$\otimes$}
\esegment
\move(1.5 -3)
\bsegment
\move(0 0)\rlvec(0.8 0)\rlvec(0 0.8)\rlvec(-0.8 0)\rlvec(0 -0.8)
\htext(0.38 0.38){$3$}
\move(1.4 0.4)\ravec(1.5 0)
 \move (9.4 -0.5)
 \clvec (9.4 -1.5)(0 -1.5)(0 -0.5)
\move(0 -0.55)\ravec(-0.1 0.2)
\htext(4.56 1.35){$0$}
\htext(2 0.7){$1$}
\esegment
\esegment

\move(5 -6)
\bsegment
\move(0 -3)
\bsegment
\move(0 0)\rlvec(0.8 0)\rlvec(0 0.8)\rlvec(-0.8 0)\rlvec(0 -0.8)
\htext(0.38 0.38){$2$}
\htext(1.15 0.4){$\otimes$}
\esegment
\move(1.5 -3)
\bsegment
\move(0 0)\rlvec(0.8 0)\rlvec(0 0.8)\rlvec(-0.8 0)\rlvec(0 -0.8)
\htext(0.38 0.38){$3$}
\esegment
\esegment

\move(10 -6)
\bsegment
\move(0 -3)
\bsegment
\move(0 0)\rlvec(0.8 0)\rlvec(0 0.8)\rlvec(-0.8 0)\rlvec(0 -0.8)
\htext(0.38 0.38){$3$}
\htext(1.15 0.4){$\otimes$}
\esegment
\move(1.5 -3)
\bsegment
\move(0 0)\rlvec(0.8 0)\rlvec(0 0.8)\rlvec(-0.8 0)\rlvec(0 -0.8)
\htext(0.38 0.38){$3$}
\esegment
\esegment
\esegment
\end{texdraw}%
\end{center}
Then
$$H(\
\raisebox{-0.3em}{\begin{texdraw}%
\fontsize{9}{9}\selectfont \drawdim em \setunitscale 1.9 \move(0
0)\rlvec(0.8 0)\rlvec(0 0.8)\rlvec(-0.8 0)\rlvec(0 -0.8)
\htext(1.1 0.2){$\otimes$} \move(1.8 0)\rlvec(0.8 0)\rlvec(0
0.8)\rlvec(-0.8 0)\rlvec(0 -0.8) \htext(0.3 0.3){$a$} \htext(2.1
0.3){$b$}
\end{texdraw}%
}\
) = \begin{cases} 1 & \qquad \hbox{\rm if} \quad a \geq b \\
0 & \qquad \hbox{\rm if} \quad a < b. \end{cases}$$

 \bi

\begin{thm}\label{thm:wtchar}\cite{(KMN)$^2$a}   Assume $\lambda \in \bar P^+$ and  $\p = (p_k)_{k=0}^\infty \in \mathcal P(\lambda)$.
Then the weight of $\p$ and the character of the irreducible $U_q(\widehat \g)$-module
$L(\lambda)$ are given by the following expressions:
\begin{eqnarray}&& \label{firsteq}  \\
&&\vspace{-.8cm} {\hbox{\rm wt}\p} = \lambda + \sum_{k=0}^\infty \left(\ov{\hbox{\rm wt}p}_k -
 \ov{\hbox{\rm wt}b}_k\right)  \nonumber  \\
&& \hspace{1.6cm}  -\left( \sum_{k=0}^\infty (k+1)\Big(H(p_{k+1} \ot p_k) - H(b_{k+1}\ot b_k)\Big)\right) \delta,  \nonumber \\
&&\hbox{\rm ch}L(\lambda) = \sum_{\p \in \mathcal P(\lambda)}  e^{\hbox{\rm wt$\p$}}.  \nonumber
\end{eqnarray}
\end{thm}
\bi
({\it Note that in Equation  (\ref{firsteq}), we are viewing $\ov{\hbox{\rm wt}p}_k$ and
 $\ov{\hbox{\rm wt}b}_k$ as classical weights, i.e. elements of the $\mathbb Z$-submodule $\mathbb Z \Lambda_0 \oplus \mathbb Z \Lambda_1 \oplus \cdots \oplus \mathbb Z \Lambda_n$ of $P$
rather than considering their restriction to $\bar Q^\vee$.})
 \bi

Since perfect crystals reveal much about the structure of crystal bases for irreducible
modules, which in turn can be used to compute their weights and characters, our goal
in the subsequent sections will be to construct perfect crystals
for all affine Lie algebras and to calculate the corresponding energy functions.

\bi
\begin{section} {A Uniform Construction of  Level 1 Perfect Crystals} \end{section}

Let  $\widehat \g$ be an affine Lie algebra and let
\begin{equation}\theta =  \begin{cases}  d_1 \alpha_1 + \cdots + d_n \alpha_n & \qquad \hbox{\rm if} \ \
\widehat \g \neq \hbox{\rm A}_{2n}^{(2)}  \\
\hf \left (d_1 \alpha_1 + \cdots +  d_n \alpha_n\right)  & \qquad \hbox{\rm if} \ \
\widehat \g = \hbox{\rm A}_{2n}^{(2)},\end{cases} \end{equation}
where the $d_i$ are as in \eqref{eq:centnull}.  Thus, when $\widehat \g = \hbox{\rm X}_n^{(1)}$
(the so-called untwisted case),  $\theta$ is the highest root of $\g$.   In all other cases except
for $\hbox{\rm A}_{2n}^{(2)}$, \  $\theta$
is the highest short root of $\g$.     The specific values
of the $d_i$ can be read from \cite[Ex.~10.1.1]{HK} or
they can be seen from the marks above the roots $\alpha_1,\dots,\alpha_n$
in  \cite[Tables Aff 1-3]{K}.  \m

Let $\mathcal B(\theta)$ denote the crystal graph of the
irreducible $U_q(\g)$-module $L_\g(\theta)$.  Thus,  the crystal graph $\mathcal
B(\theta)$ corresponds to the adjoint representation of
$\g$ in the untwisted case and to the ``little'' adjoint
representation of $\g$ (with highest weight the highest short
root) in all other cases but $\hbox{\rm A}_{2n}^{(2)}$. The set 
$\Lambda := \hbox{\rm wt}\,\mathcal B(\theta)$ of weights of
$\mathcal B(\theta)$ is a subset of $\Phi \cup \{0\}$, where
$\Phi$ is the root system of $\g$ (except when $\widehat \g$ is of
type A$_{2n}^{(2)}$). In the untwisted case  equality holds,
$\Lambda = \Phi \cup \{0\}$. \m

Let $\Phi^+$ and $\Phi^- = -\Phi^+$ denote the positive and
negative roots respectively of $\g$. Set $\Lambda^+ = \Lambda \cap
\Phi^+$, \ $\Lambda^- = - \Lambda^+$, so that $\Lambda =
\Lambda^+\cup \{0\} \cup \Lambda^-$.   Note if  $\widehat \g = \hbox{\rm
A}_{2n}^{(2)}$,  then 
$$\Lambda = \Lambda^+ \cup \Lambda^-
=\left \{ \pm( \alpha_i + \cdots + \alpha_{n-1} + {\textstyle \hf} \alpha_n) \ \big | \ 
i=1, \dots, n-1\right \}\cup \left \{ \pm  {\textstyle \hf} \alpha_n\right \}.$$ \m

Correspondingly, we write

\begin{equation} \mathcal B(\theta) =
\begin{cases} \{x_\alpha \mid \alpha \in \Lambda^+ \} \cup \{y_i \mid
\alpha_i \in \Lambda^+\} \cup \{x_{-\alpha} \mid \alpha \in
\Lambda^+\} & \quad \text{if} \ \ \widehat \g \neq \hbox{\rm
A}_{2n}^{(2)}, \\
\{x_\alpha \mid \alpha \in \Lambda^+ \} \cup \{x_{-\alpha} \mid
\alpha \in \Lambda^+\} & \quad \text{if} \ \ \widehat \g =
\hbox{\rm A}_{2n}^{(2)}.
\end{cases}
\end{equation}

\noindent Hence in the untwisted case, $\mathcal B(\theta) = \{x_{\pm \alpha} \mid \alpha \in \Phi^+\}
\cup \{y_i \mid i=1,\dots,n\}$.

Set $\mathcal B(0) = \{\emptyset\}$, which we identify with  the
crystal graph of the one-dimensional $U_q(\g)$-module $L_\g(0)$. As we
argue below, the set

\begin{equation} \mathcal B = \mathcal B(\theta) \sqcup \mathcal B(0) \end{equation}

\noindent can be endowed with  a crystal structure as follows:

\begin{equation}\label{eq:crys} \begin{array}{cccc}
&(i \neq 0)   &\qquad x_\alpha\ \iarrow \ x_\beta \ \Longleftrightarrow \ \alpha-\alpha_i = \beta & \qquad (\alpha,\beta \in
\Lambda), \\
& &\qquad x_{\alpha_i} \ \iarrow \ y_i \ \iarrow \ x_{-\alpha_i} & \qquad (\alpha_i \in \Lambda^+),\\
&(i=0)  &\qquad x_\alpha \ \zarrow \ x_{\beta} \ \Longleftrightarrow \ \alpha+\theta = \beta &
\qquad (\alpha,\beta \neq \pm \theta), \\
& &\qquad x_{-\theta} \ \zarrow \ \emptyset \ \zarrow \ x_{\theta}. & \end{array}\end{equation}

We remark that in the $\widehat \g = \hbox{\rm A}_{2n}^{(2)}$ case, no $\alpha_i$ belongs to $\Lambda^+$.
Now we are ready to state our main theorem.
\bi
\begin{thm}\label{thm:main}   $\mathcal B = \mathcal B(\theta) \sqcup \mathcal B(0)$
with the structure given in  \eqref{eq:crys}  is a perfect
crystal of level 1 for every quantum affine algebra $U_q'(\widehat \g)$.  \end{thm}
\bi

Before embarking on the proof of Theorem \ref{thm:main},  we present several examples.
In doing so, we use $\varpi_1,\varpi_2, \dots$ to denote the fundamental weights of
the finite-dimensional algebra $\g$.
\bi
\begin{exs}

{\rm (1)} \ $\widehat \g = \hbox{\rm A}_{2}^{(1)}$, $\g = \hbox{\rm A}_2$, $\theta = \alpha_1+\alpha_2 = \varpi_1+\varpi_2$, \hbox{\rm and} $c = h_0+h_1+h_2$.
\bi
\begin{center}
\begin{texdraw}%
\drawdim in
\arrowheadsize l:0.065 w:0.03
\arrowheadtype t:F
\fontsize{11}{11}\selectfont
\textref h:C v:C
\drawdim em
\setunitscale 1.9
\htext(0 0){$x_{\theta}$}
\htext(2 0){$x_{\alpha_2}$}
\htext(4 0){$y_2$}
\htext(0 -1.5){$x_{\alpha_1}$}
\htext(2 -1.5){$y_1$}
\htext(4 -1.5){$x_{-\alpha_2}$}
\htext(2 -3){$x_{-\alpha_1}$}
\htext(4 -3){$x_{-\theta}$}
\htext(0 -3){$\emptyset$}
\move(0.4 0)\ravec(1.1 0)
\move(2.4 0)\ravec(1.1 0)
\move(0.4 -1.5)\ravec(1.1 0)
\move(2.4 -3)\ravec(1.1 0)
\move(0 -0.4)\ravec(0 -0.9)
\move(4 -0.4)\ravec(0 -0.9)
\move(2 -1.9)\ravec(0 -0.9)
\move(4 -1.9)\ravec(0 -0.9)
 \move (2.3 -2.8)
 \clvec (2.6 -2.8)(2.6 -0.3)(2.3 -0.3)
 \move (4 -1.8)
 \clvec (4 -2.1)(0 -2.1)(0 -1.8)
 \move(0 -1.8)\ravec(-0.1 0.08)
 \move(2.32 -0.3)\ravec(-0.06 0.1)
 \move (4 -3.4)
 \clvec (4 -3.7)(0.2 -3.7)(0.2 -3.4)
 \move (-0.4 -3)
 \clvec (-0.7 -3)(-0.7 -0.2)(-0.4 -0.2)
 \move(0.21 -3.42)\ravec(-0.2 0.2)
 \move(-0.4 -0.2)\ravec(0.1 0.1)
\htext(0.8 0.3){\tiny$1$}
\htext(2.8 0.3){\tiny$2$}
\htext(0.8 -1.2){\tiny$1$}
\htext(2.8 -2.7){\tiny$2$}
\htext(-0.3 -0.75){\tiny$2$}
\htext(3.7 -0.75){\tiny$2$}
\htext(3.7 -2.3){\tiny$1$}
\htext(1.7 -2.3){\tiny$1$}
\htext(-0.9 -1.5){\tiny$0$}
\htext(2.8 -1.5){\tiny$0$}
\htext(0.8 -2.3){\tiny$0$}
\htext(0.8 -3.9){\tiny$0$}
\htext(9 0){$b^{\Lambda_0}=b_{\Lambda_0}=\emptyset$}
\htext(9.15 -1){$b^{\Lambda_1}=b_{\Lambda_1}=y_1$}
\htext(9.17 -2){$b^{\Lambda_2}=b_{\Lambda_2}=y_2$}
\end{texdraw}%
\end{center}
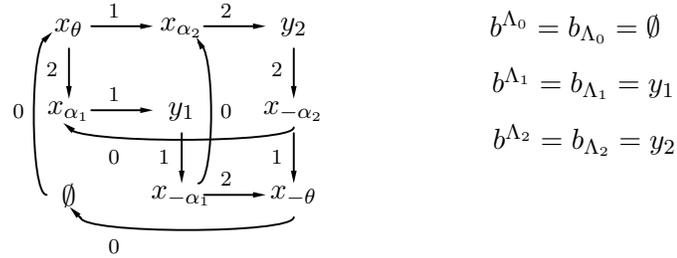
\bi
{\rm (2)} $\widehat \g = \hbox{\rm D}_{4}^{(3)}$, $\g = \hbox{\rm G}_2$, $\theta = 2\alpha_1+\alpha_2 = \varpi_1$, \hbox{\rm and} $c = h_0+2h_1+3h_2$.
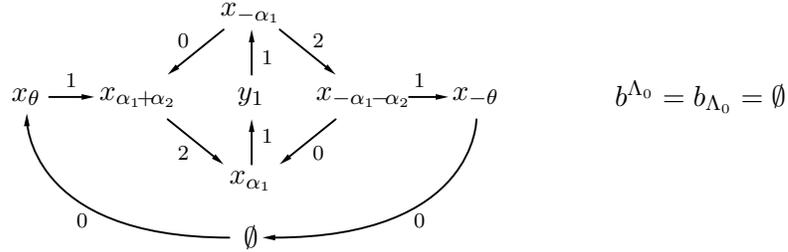
\bi
\begin{center}
\begin{texdraw}%
\drawdim in
\arrowheadsize l:0.065 w:0.03
\arrowheadtype t:F
\fontsize{11}{11}\selectfont
\textref h:C v:C
\drawdim em
\setunitscale 1.9
\htext(0 1.5){$x_{-\alpha_1}$}
\htext(-2 0){$x_{\alpha_1\!+\!\alpha_2}$}
\htext(0 0){$y_1$}
\htext(2 0){$x_{-\alpha_1\!-\!\alpha_2}$}
\htext(4 0){$x_{-\theta}$}
\htext(-4 0){$x_{\theta}$}
\htext(0 -1.5){$x_{\alpha_1}$}
\htext(0 -2.5){$\emptyset$}
\move(0 0.4)\ravec(0 0.8)
\move(0 -1.2)\ravec(0 0.8)
\move(-3.6 0)\ravec(0.8 0)
\move(2.8 0)\ravec(0.7 0)
\move(-1.5 -0.4)\ravec(1 -0.8)
\move(1.5 -0.4)\ravec(-1 -0.8)
\move(-0.5 1.2)\ravec(-1 -0.8)
\move(0.5 1.2)\ravec(1 -0.8)
 \move (-0.4 -2.5)
 \clvec (-4 -2.5)(-4 -0.4)(-4 -0.4)
 \move (0.4 -2.5)
 \clvec (4 -2.5)(4 -0.4)(4 -0.4)
 \move(-4 -0.5)\ravec(0 0.2)
 \move(0.4 -2.5)\ravec(-0.2 0)
\htext(3 -2.2){\tiny$0$}
\htext(-3 -2.2){\tiny$0$}
\htext(-3.2 0.3){\tiny$1$}
\htext(3 0.3){\tiny$1$}
\htext(0.3 0.7){\tiny$1$}
\htext(0.3 -0.7){\tiny$1$}
\htext(1.2 1){\tiny$2$}
\htext(-1.2 1){\tiny$0$}
\htext(1.2 -1){\tiny$0$}
\htext(-1.2 -1){\tiny$2$}
\htext(8 0){$b^{\Lambda_0}=b_{\Lambda_0}=\emptyset$}
\end{texdraw}%
\end{center}
\bi
{\rm (3)} $\widehat \g = \hbox{\rm C}_{2}^{(1)}$, $\g = \hbox{\rm C}_2$, $\theta = 2\alpha_1+\alpha_2 = 2\varpi_1$, \hbox{\rm and} $c = h_0+h_1+h_2$.
\bi
\begin{center}
\begin{texdraw}%
\drawdim in
\arrowheadsize l:0.065 w:0.03
\arrowheadtype t:F
\fontsize{11}{11}\selectfont
\textref h:C v:C
\drawdim em
\setunitscale 1.9
\htext(0 0){$x_{\theta}$}
\htext(2 0){$x_{\alpha_1\!+\!\alpha_2}$}
\htext(4 0){$x_{\alpha_2}$}
\htext(2 -1.5){$x_{\alpha_1}$}
\htext(3.3 -1.5){$y_1$}
\htext(4 -1.1){$y_2$}
\htext(6 -1.5){$x_{-\alpha_1}$}
\htext(4 -3){$x_{-\alpha_2}$}
\htext(6 -3){$x_{-\alpha_1\!-\!\alpha_2}$}
\htext(8 -3){$x_{-\theta}$}
\htext(2 -4.5){$\emptyset$}
\move(0.4 0)\ravec(0.8 0)
\move(2.7 0)\ravec(0.8 0)
\move(2.4 -1.5)\ravec(0.6 0)
\move(3.7 -1.5)\ravec(1.6 0)
\move(4.6 -3)\ravec(0.5 0)
\move(6.9 -3)\ravec(0.6 0)
\move(2 -0.4)\ravec(0 -0.9)
\move(4 -0.4)\ravec(0 -0.5)
\move(4 -1.4)\ravec(0 -1.4)
\move(6 -1.9)\ravec(0 -0.9)
 \move (2 0.3)
 \clvec (3 1)(5.5 1)(6 -1.2)
 \move (6 -3.3)
 \clvec (5 -4)(2.5 -4)(2 -1.8)
 \move(2 0.3)\ravec(-0.1 -0.08)
 \move(2 -1.8)\ravec(-0.02 0.1)
 \move (2.5 -4.5)
 \clvec (8 -4.5)(8 -4.5)(8 -3.4)
 \move (-0 -0.4)
 \clvec (-0 -4.5)(0 -4.5)(1.6 -4.5)
 \move(2.5 -4.5)\ravec(-0.2 0)
 \move(0 -0.4)\ravec(0 0.1)
\htext(5 -4.8){\tiny$0$}
\htext(-0.3 -2){\tiny$0$}
\htext(5.2 0.6){\tiny$0$}
\htext(2.3 -3){\tiny$0$}
\htext(0.8 0.3){\tiny$1$}
\htext(3 0.3){\tiny$1$}
\htext(2.6 -1.2){\tiny$1$}
\htext(4.6 -1.2){\tiny$1$}
\htext(4.7 -2.7){\tiny$1$}
\htext(7 -2.7){\tiny$1$}
\htext(1.7 -0.7){\tiny$2$}
\htext(3.7 -0.7){\tiny$2$}
\htext(3.7 -2.2){\tiny$2$}
\htext(5.7 -2.2){\tiny$2$}
\htext(11 0){$b^{\Lambda_0}=b_{\Lambda_0}=\emptyset$}
\htext(11.15 -1){$b^{\Lambda_1}=b_{\Lambda_1}=y_1$}
\htext(11.17 -2){$b^{\Lambda_2}=b_{\Lambda_2}=y_2$}
\end{texdraw}%
\end{center}
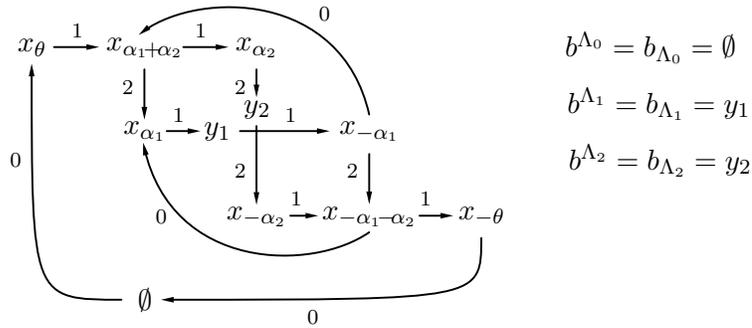
\bi
{\rm (4)} $\widehat \g = \hbox{\rm A}_{4}^{(2)}$, $\g = \hbox{\rm C}_2$, $\theta = \alpha_1+\hf \alpha_2 = \varpi_1$, \hbox{\rm and} $c = h_0+2h_1+2h_2$.
\bi
\begin{center}
\begin{texdraw}%
\drawdim in
\arrowheadsize l:0.065 w:0.03
\arrowheadtype t:F
\fontsize{11}{11}\selectfont
\textref h:C v:C
\drawdim em
\setunitscale 1.9
\htext(0 0){$x_{\theta}$}
\htext(2 0){$x_{\hf\alpha_2}$}
\htext(4.2 0){$x_{-\hf\alpha_2}$}
\htext(7 0){$x_{-\alpha_1\!-\!\hf\alpha_2}$}
\htext(3.5 -1.5){$\emptyset$}
\move(0.4 0)\ravec(1 0)
\move(2.6 0)\ravec(0.8 0)
\move(5 0)\ravec(1 0)
 \move (7 -0.4)
 \clvec (7 -1.5)(4 -1.5)(4 -1.5)
 \move (0 -0.4)
 \clvec (0 -1.5)(3 -1.5)(3 -1.5)
 \move(4 -1.5)\ravec(-0.2 0)
 \move(-0.01 -0.4)\ravec(-0.02 0.1)
\htext(5.5 -1.7){\tiny$0$}
\htext(1.5 -1.7){\tiny$0$}
\htext(0.8 0.3){\tiny$1$}
\htext(2.9 0.3){\tiny$2$}
\htext(5.4 0.3){\tiny$1$}
\htext(11 0){$b^{\Lambda_0}=b_{\Lambda_0}=\emptyset$}
\end{texdraw}%
\end{center}
\bi
{\rm More generally, for arbitrary $n$ we have}\m

{\rm (5)} $\widehat \g = \hbox{\rm A}_{2n}^{(2)}$, $\g = \hbox{\rm
C}_n$, $\theta = \alpha_1+\alpha_2 + \cdots + \alpha_{n-1}+\hf
\alpha_n = \varpi_1$, \hbox{\rm and} $c = h_0+2h_1+2h_2 + \cdots+2h_n$. \bi
\begin{center}
\begin{texdraw}%
\drawdim in
\arrowheadsize l:0.065 w:0.03
\arrowheadtype t:F
\fontsize{11}{11}\selectfont
\textref h:C v:C
\drawdim em
\setunitscale 1.9
\htext(-7 0){$x_{\theta}$}
\htext(-5 0){$x_{\theta\!-\!\alpha_1}$}
\htext(-3 0){$\cdots$}
\htext(-1 0){$x_{\hf\alpha_n}$}
\htext(1 0){$x_{-\hf\alpha_n}$}
\htext(3 0){$\cdots$}
\htext(5 0){$x_{-\theta\!+\!\alpha_1}$}
\htext(7 0){$x_{-\theta}$}
\htext(0 -1.5){$\emptyset$}
 \move (-0.4 -1.5)
 \clvec (-7 -1.5)(-7 -0.5)(-7 -0.5)
 \move (0.4 -1.5)
 \clvec (7 -1.5)(7 -0.5)(7 -0.5)
 \move(0.5 -1.5)\ravec(-0.2 0)
 \move(-7 -0.5)\ravec(-0.05 0.1)
\move(-6.6 0)\ravec(0.9 0)
\move(-4.4 0)\ravec(0.8 0)
\move(-2.3 0)\ravec(0.7 0)
\move(-0.4 0)\ravec(0.6 0)
\move(1.7 0)\ravec(0.7 0)
\move(3.6 0)\ravec(0.6 0)
\move(5.7 0)\ravec(0.9 0)
\htext(4 -1.8){\tiny$0$}
\htext(-4 -1.8){\tiny$0$}
\htext(-6.3 0.3){\tiny$1$}
\htext(-4.2 0.3){\tiny$2$}
\htext(-2 0.3){\tiny$n\!\!-\!\!1$}
\htext(-0.2 0.3){\tiny$n$}
\htext(2 0.3){\tiny$n\!\!-\!\!1$}
\htext(6.1 0.3){\tiny$1$}
\htext(3.9 0.3){\tiny$2$}
\htext(9.5 0){$b^{\Lambda_0}=b_{\Lambda_0}=\emptyset$}
\end{texdraw}%
\end{center}

\bi \noindent {\rm

}
\end{exs} \bi

{\it Proof of Theorem \ref {thm:main}}.  \   Our argument  will proceed in a series of steps in which we verify
that $\mathcal B$ satisfies the conditions in Definition \ref{defn:pc}.

\bi
\noindent \textbf {Step 1.}  \  {\it The space
$$V = \left (\bigoplus_{\alpha \in \Lambda^+} \mathbb C(q) x_{\alpha}
\right)
\oplus \Bigg(\bigoplus_{i, \ {\alpha_i \in \Lambda^+}} \mathbb C(q) y_i
\Bigg)
\oplus \left (\bigoplus_{\alpha \in \Lambda^+} \mathbb C(q) x_{-\alpha}
\right) \oplus \mathbb C(q)\, \emptyset$$
can be given the structure of a $U_q'(\widehat \g)$-module. }
\m

To show this,   we must assign an action of the generators of $U_q'(\widehat \g)$ on $V$:

\begin{eqnarray*}
&& q^h.x_\beta = q^{\beta(h)}x_\beta , \qquad q^h y_i = y_i, \qquad
q^h.\emptyset = \emptyset; \\
&& e_i.x_\beta =
\begin{cases}[\varphi_i(x_\beta)+1]_{q_i} x_{\beta+\alpha_i} & \qquad \hbox{\rm if} \ \beta+\alpha_i \in \Lambda \\
0 & \qquad \hbox{\rm otherwise}  \end{cases}  \qquad   (i \neq 0);  \\
&&f_i.x_\beta = \begin{cases}[\varepsilon_i(x_\beta)+1]_{q_i} x_{\beta-\alpha_i} & \qquad \hbox{\rm if} \ \beta-\alpha_i \in \Lambda \\
0 & \qquad \hbox{\rm otherwise}  \end{cases}  \qquad  (i \neq 0);  \\
&&e_i.x_{-\alpha_i} = y_i , \qquad e_i.y_i = [2]_{q_i} x_{\alpha_i}  \qquad \  (i \neq 0);  \\
&&f_i.x_{\alpha_i} = y_i ,\qquad \ f_i.y_i = [2]_{q_i} x_{-\alpha_i}  \qquad  (i \neq 0);  \\
&&e_0.x_{\beta} = \begin{cases} x_{\beta-\theta} & \qquad \hbox{\rm if} \ \beta-\theta  \in \Lambda \\
0 & \qquad \hbox{\rm otherwise} \end{cases} \\
&&f_0.x_{\beta} = \begin{cases} x_{\beta+\theta} & \qquad \hbox{\rm if} \ \beta+\theta  \in \Lambda \\
0 & \qquad \hbox{\rm otherwise} \end{cases} \\
&&e_0.x_\theta = \emptyset, \qquad \,e_0.\emptyset = [2]_{q_0} \, x_{-\theta}, \\
&&f_0.x_\theta = \emptyset, \qquad \, f_0.\emptyset = [2]_{q_0} \, x_{\theta},
\end{eqnarray*}
for all $\beta \in \Lambda^{\pm}$.   By \cite{KMPY}, it suffices to check that relations (1)-(4)  in
Definition \ref{defn:qalg} hold.    This can be done on a case-by-case basis.
Since the action of $U_q(\g)$ on $V$ is the same as on the
$U_q(\g)$-module  $L_\g(\theta) \oplus L_\g(0)$, the relations that need to be verified are
the ones involving $e_0,f_0$, and $K_0^{\pm 1}$.
Here are a few sample calculations.   In the second one, we will use
the fact that $\beta(h_0) \in \{-1,0,1\}$ for all $\beta \in \Lambda$ when
$\widehat \g \neq \hbox{\rm A}_{2n}^{(2)}$.
If $\beta(h_0) = -1$, then $\beta -\theta \in \Lambda$  but $\beta+\theta \not \in \Lambda$,
and analogously, if $\beta(h_0) = 1$, then $\beta +\theta \in \Lambda$  but $\beta - \theta
\not \in \Lambda$.

\begin{eqnarray*} &&(e_0f_0-f_0e_0).\emptyset = [2]_{q_0} e_0.x_{\theta}
-[2]_{q_0} f_0.x_{-\theta} =  [2]_{q_0} \bigl(\emptyset -\emptyset) = 0, \\
&&\left(\frac{K_0 - K_0^{-1}}{q_0-q_0^{-1}}\right).\emptyset = 0,  \\
&&(e_0f_0-f_0e_0).x_{\beta} =  \begin{cases}
-f_0e_0.x_\beta = -x_\beta &\qquad \hbox{\rm if} \ \ \beta(h_0) = -1, \\
e_0f_0.x_\beta = x_\beta &\qquad \hbox{\rm if} \ \ \beta(h_0) = 1, \\
0 &\qquad \hbox{\rm otherwise}, \end{cases}\\
&&\left(\frac{K_0 - K_0^{-1}}{q_0-q_0^{-1}}\right).x_{\beta} = \frac{q_0^{\beta(h_0)} - q_0^{-\beta(h_0)}}
{q_0-q_0^{-1}}x_\beta  =  \begin{cases}
-x_\beta &\qquad \hbox{\rm if} \ \ \beta(h_0) = -1, \\
x_\beta &\qquad \hbox{\rm if} \ \ \beta(h_0) = 1, \\
0 &\qquad \hbox{\rm otherwise}. \end{cases}\\
 \end{eqnarray*}

\bi
\noindent \textbf{Step 2.} \  {\it There exists a classical weight $\lambda_0 \in
\bar P$ such that $|\, \mathcal B_{\lambda_0}\,| = 1$ and $\hbox{\rm wt}(\mathcal B)
\subset \lambda_0+ \frac{1}{d_0} \sum_{i \neq 0} \mathbb Z_{\leq 0} \alpha_i$.}  \m

This is easily seen by taking $\lambda_0 = \theta$.

\bi
\noindent \textbf{Step 3.} \  {\it  For all $b \in \mathcal B$, we have $\la c, \varepsilon(b) \ra \geq 1$.}
\m

First suppose that $b \in \mathcal B(\theta)$ and $b \neq
x_\theta$.  Then there exists an $i \neq 0$ such that $\tilde e_i
b \neq 0$;  i.e.,  $\varepsilon_i(b) \geq 1$, so that $\la c,
\varepsilon(b) \ra \geq 1$ must hold. When $b = x_\theta$,  we
have $\varepsilon_0(b) = 2$  (see \eqref{eq:crys}) so that $\la c,
\varepsilon(b) \ra \geq 1$.   Finally, when $b = \emptyset$, then
$\varepsilon_0(\emptyset) = 1$, so that $\la c, \varepsilon(b) \ra
\geq 1$ in this case also.

\bi
\noindent \textbf{Step 4.} \  {\it   For all $\lambda \in P^+$ with $\lambda(c) = 1$, there exist
unique elements $b^\lambda$ and $b_\lambda \in \mathcal B$ such that
$\varepsilon(b^\lambda) = \lambda$ and $\varphi(b_\lambda) = \lambda$.}
\m

When $\lambda = \Lambda_0$, we can take $b^{\Lambda_0} = \emptyset = b_{\Lambda_0}$.
Now when $\lambda = \Lambda_i$ for $i \neq 0$, then setting
$b^{\Lambda_i} = y_i = b_{\Lambda_i}$ will give the desired result.

\bi
All that remains in the proof of Theorem \ref{thm:main}
is to show that the crystal graph  $\mathcal B \ot \mathcal B$ is connected.
We devote the next section to this task.

\bi
\begin{section} {$\mathcal B \ot \mathcal B$ is connected} \end{section}

Our approach to proving this can be summarized as follows.  We forget the $0$-arrows in $\mathcal B \ot \mathcal B$
and view it as a crystal graph for the quantum algebra $U_q(\g)$ associated to the
simple Lie algebra $\g$:
\begin{equation}\label{eq:dec}\mathcal B \ot \mathcal B = \big(\mathcal B(\theta) \ot \mathcal B(\theta)\big)
\sqcup \big(\mathcal B(\theta) \ot \mathcal B(0)\big) \sqcup \big(\mathcal B(0) \ot \mathcal B(\theta)\big)
\sqcup \big(\mathcal B(0) \ot \mathcal B(0)\big).\end{equation}
Since crystals corresponding to simple modules are connected, it suffices to
locate the maximal vectors  ($\tilde e_i b = 0$ for all $i \in \mathcal I \setminus \{0\}$)
inside the components on the right and show that they are all connected to one another
by various $i$-arrows for $i \in \mathcal I$.  \m

There are obvious maximal vectors inside $\mathcal B \ot \mathcal B$,

\begin{itemize}
\item[{\rm (1)}]  $x_\theta \ot x_\theta$
\item[{\rm (2)}]   $x_\theta \ot \emptyset$
\item[{\rm (3)}]   $\emptyset \ot x_\theta$
\item[{\rm (4)}]   $\emptyset \ot \emptyset$
\item[{\rm (5)}]  $x_\theta \ot x_{-\theta}$,
\end{itemize}
and they can be connected as displayed below:

$$
\begin{array}{ccccccccccc} \emptyset \ot \emptyset &\zarrow&  \emptyset \ot x_\theta &\dotarrow& \emptyset \ot x_{-\theta} &\zarrow& x_\theta \ot x_{-\theta} &\zarrow& x_\theta \ot \emptyset &\zarrow& x_\theta \ot x_\theta,  \end{array}$$
where $\dotarrow$ indicates that an appropriate sequence of Kashiwara operators  $\tilde f_i$ with  $i \in \mathcal I \setminus  \{0\}$ has been applied.   All other maximal vectors have the form
\begin{itemize}
\item[{\rm (6)}]  $x_\theta \ot x_{\theta-\alpha}$   for some $\alpha \in \Lambda^+$ or
\item[{\rm (7)}]   $x_\theta \ot y_i$  for some $i$ such that $\alpha_i \in
\Lambda^+$. (When $\widehat \g$ is of type $\hbox{\rm A}_{2n}^{(2)}$ or $\hbox {\rm D}_{n+1}^{(2)}$, this
case does not occur,  because no $\alpha_i$ belongs to $\Lambda^+$.)
\end{itemize}
Since they can be connected to $\emptyset \ot x_{\theta}$ via
$$
\begin{array}{ccccccccccc}x_\theta \ot  x_{\theta-\alpha} &\zlarrow& x_\theta  \ot x_{-\alpha} &\zlarrow& \emptyset \ot x_{-\alpha} &\dotlarrow& \emptyset  \ot x_{\theta} &\dotarrow& \emptyset  \ot y_i &\zarrow& x_\theta \ot y_i,  \end{array}$$
the entire crystal graph $\mathcal B \ot \mathcal B$ must be connected.   \qed
\bi

\begin{section}{Crystal Algebras} \end{section}

The irreducible $U_q(\g)$-module $L_\g(\theta)$ is special,  since except
for type A$_n$, $n \geq 2$,  there is a unique (up to scalar factor)
$U_q(\g)$-invariant projection,    

\begin{equation}\label{eq:mult}  \mathfrak m:  L_\g(\theta) \otimes L_\g(\theta)
\rightarrow L_\g(\theta),   \end{equation}
 
 \noindent  giving  a canonical algebra structure  on $L_\g(\theta)$.   
(Imposing skew-symmetry will also ensure that such a map is unique for type A$_n$, $n \geq 2$.)
When $\theta$ is the highest root of $\g$ and  $q=1$ (i.e. in the classical  $\g$-module case),  the module $L_q(\g)$  is  isomorphic to the
adjoint module $\g$  and \eqref{eq:mult} is just the  Lie algebra structure.  
When $\theta$ is the highest short root of $\g$ (and again $q=1$),
we obtain the 27-dimensional exceptional simple Jordan algebra 
minus its unit element  if $\g$ is of type F$_4$,  and the
8-dimensional octonion algebra minus its unit element  if $\g$ is of type G$_2$.      Thus,
\eqref{eq:mult} defines  quantum
analogues of the above algebras for generic values of $q$.   Taking the limit
$q=0$,  we obtain a strict morphism of crystals: 

\begin{equation}\label{eq:morph}  \mathfrak m:  \mathcal B(\theta) \ot \mathcal B(\theta)
\rightarrow \mathcal B(\theta).  \end{equation}

\noindent  We will describe this morphism (and its inverse)  explicitly in terms
of  the root vectors, since
it is of independent interest.   It will also be crucial in
determining  the energy function $H$.

Our computations of the energy function $H$ in the next section  will rely heavily  on knowing the
various connected components that result from omitting the 0-arrows of $\mathcal B \otimes \mathcal B$, since $H$ must be  constant on those components. 
The  connected components  $\mathcal C(x_\theta \ot y_i)$, for values of  $i$ 
such that  vertex $i$ is connected to vertex 0 in the Dynkin diagram of $\widehat \g$,  are
of particular interest, for they are isomorphic  to the crystal $\mathcal B(\theta)$.
In the next proposition, we construct the isomorphism $\Psi: \mathcal B(\theta) \rightarrow 
\mathcal C(x_\theta \ot y_i)$  explicitly.  
The strict morphism  $\mathfrak m:  \mathcal B(\theta) \ot \mathcal B(\theta)
\rightarrow \mathcal B(\theta)$  is the inverse  of $\Psi$ 
on  $\mathcal C(x_\theta \ot y_i)$ and is zero on all other components of
$\mathcal B(\theta) \ot \mathcal B(\theta)$.   
We regard the morphism $\mathfrak m$  as a multiplication on $\mathcal B(\theta)$
resulting in a ``crystal algebra''   $(\mathcal B(\theta), \mathfrak m)$.   We will discuss this
algebra and provide an interesting example.

 \bi
\begin{prop}\label{thetacomp} There is a crystal isomorphism $\mathcal B(\theta) \iso \mathcal C(x_\theta \ot y_i)$
whenever  vertex $i$ is connected to vertex 0 in the Dynkin diagram of $\widehat \g$
and $\alpha_i$ belongs to $\Lambda^+$.
\end{prop}

\bi
\begin{rem}  When $\widehat \g$ is of type $\hbox{\rm A}_{2n}^{(2)}$ or $\hbox{\rm D}_{n+1}^{(2)}$,  there
is no value of $i$ satisfying the hypotheses of the proposition, and there are no
components of $\mathcal B(\theta) \ot \mathcal B(\theta)$ isomorphic to $\mathcal B(\theta)$.  \end{rem}
\bi

\begin{proof}
Our proof of Proposition \ref{thetacomp}  will be broken into a number of cases.
We first establish that a crystal isomorphism  $\mathcal B(\theta)
\iso \mathcal C(x_\theta \ot y_i)$ exists for   $\widehat \g \neq
\hbox{\rm A}_n^{(1)}$, $\hbox{\rm C}_n^{(1)}$.
For such $\widehat \g$,  the root  $\theta$ is equal to
the fundamental weight $\varpi_i$.   We begin the argument by introducing a ``grading''
on $\Lambda^+$ as follows:   set
$$\Lambda_j^+  = \{ \gamma \in \Lambda^+ \mid
\hbox{\rm coefficient  of} \ \alpha_i \ \hbox{\rm in} \ \gamma \ \hbox{\rm is} \ j\}.$$

\noindent  Then $\Lambda^+ = \Lambda_0^+ \sqcup \Lambda_1^+ \sqcup \Lambda_2^+$ where $\Lambda_2^+ = \{ \theta\}$.      We construct a map
$\Psi:  \mathcal B(\theta) \rightarrow  \mathcal C(x_\theta \ot y_i)$ using this decomposition.
First,  for  $\gamma = \theta \in \Lambda_2^+$,  set 

\begin{equation}\label{eq:first} \Psi(x_\gamma) = x_\theta \ot y_i.  \end{equation}

Next, if $\gamma \in \Lambda_1^+$,  then there is a sequence $(i_1 = i, i_2, \dots, i_\ell)$ with
$1 \leq i_\ell \neq i \leq n$ for $\ell > 1$ so  that $\tilde f_{i_\ell} \cdots \tilde f_{i_2} \tilde f_{i_1} x_\theta = x_\gamma$.
Set $\beta := \alpha_{i_1} + \cdots + \alpha_{i_\ell}$.   Then by the tensor product rules,
$\tilde f_{i_\ell} \cdots \tilde f_{i_2} \tilde f_{i_1} (x_\theta \ot \emptyset) = x_\theta \ot x_{-\beta}
\in \mathcal C(x_\theta \ot y_i) \subset \mathcal B(\theta) \ot \mathcal B(\theta)$, and we 
define
\begin{equation}\label{eq:second}  \Psi(x_\gamma) =  x_\theta \ot x_{-\beta}.
\end{equation}     

Finally, suppose $\gamma \in \Lambda_0^+$,
and let  $\hbox{\rm supp}(\gamma)$ denote the set of indices $k$ such that $\alpha_k$ occurs in $\gamma$
with nonzero coefficient.  We take
$\un j = (j_1, \dots, j_t =i)$ to be  the sequence of nodes in the Dynkin diagram connecting the nodes
corresponding to the indices  in $\hbox{\rm supp}(\gamma)$
with $i$  as  in the examples below.

\begin{exs} {\rm (1)}  The root $\gamma = \alpha_1 + \alpha_2$ for $\hbox{\rm E}_6^{(1)}$ has
$\hbox{\rm supp}(\gamma) = \{1,2\}$.   Since $\theta = \varpi_6$, we have
$\un j = (3,6)$.  \m

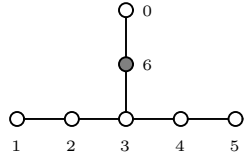
\begin{center}
\begin{texdraw}%
\drawdim in
\arrowheadsize l:0.065 w:0.03
\arrowheadtype t:F
\fontsize{6}{6}\selectfont
\textref h:C v:C
\drawdim em
\setunitscale 1.9
\move(0 0)\lcir r:0.2
\move(1.5 0)\lcir r:0.2
\move(3 0)\lcir r:0.2
\move(4.5 0)\lcir r:0.2
\move(6 0)\lcir r:0.2
\move(3 1.5)\fcir f:0.5 r:0.2 \lcir r:0.2
\move(3 3)\lcir r:0.2
\move(0.2 0)\rlvec(1.1 0)
\move(1.7 0)\rlvec(1.1 0)
\move(3.2 0)\rlvec(1.1 0)
\move(4.7 0)\rlvec(1.1 0)
\move(3 0.2)\rlvec(0 1.1)
\move(3 1.7)\rlvec(0 1.1)
\htext(0 -0.7){$1$}
\htext(1.5 -0.7){$2$}
\htext(3 -0.7){$3$}
\htext(4.5 -0.7){$4$}
\htext(6 -0.7){$5$}
\htext(3.6 1.5){$6$}
\htext(3.6 3){$0$}
\end{texdraw}%
\end{center}  \m

{\rm (2)} The root $\gamma = \alpha_2 + 2 \alpha_3$  for $\hbox{\rm F}_4^{(1)}$ has
$\hbox{\rm supp}(\gamma) = \{2,3\}$.   Since $\theta = \varpi_1$, the sequence
$\un j$ is a singleton $\un j = (1)$.  \m

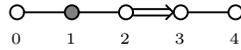
\begin{center}
\begin{texdraw}%
\drawdim in
\arrowheadsize l:0.065 w:0.03
\arrowheadtype t:F
\fontsize{6}{6}\selectfont
\textref h:C v:C
\drawdim em
\setunitscale 1.9
\move(0 0)\lcir r:0.2
\move(1.5 0)\fcir f:0.5 r:0.2 \lcir r:0.2
\move(3 0)\lcir r:0.2
\move(4.5 0)\lcir r:0.2
\move(6 0)\lcir r:0.2
\move(0.2 0)\rlvec(1.1 0)
\move(1.7 0)\rlvec(1.1 0)
\move(4.7 0)\rlvec(1.1 0)
\move(3.2 0.08)\rlvec(1 0)
\move(3.2 -0.08)\rlvec(1 0)
\move(4.3 0)\rlvec(-0.2 0.2)
\move(4.3 0)\rlvec(-0.2 -0.2)
\htext(0 -0.7){$0$}
\htext(1.5 -0.7){$1$}
\htext(3 -0.7){$2$}
\htext(4.5 -0.7){$3$}
\htext(6 -0.7){$4$}
\end{texdraw}%
\end{center}  \m
\end{exs}

Once we have the sequence $\un j$ for $\gamma \in \Lambda_0^+$, we set  $\alpha = \alpha_{j_1} + \cdots + \alpha_{j_t} = \alpha_{j_1}
+ \cdots + \alpha_{j_{t-1}} + \alpha_i$ so that $\gamma + \alpha \in \Lambda_1^+$.
Then we may write $\gamma + \alpha = \theta - \beta$.  Note $\beta\neq 0$, since
the coefficient of $\alpha_i$ in $\gamma + \alpha$ is 1.   We
set 

\begin{equation}\label{eq:third}  \Psi(x_\gamma) = x_{\theta-\alpha} \ot x_{-\beta}.
\end{equation}   This is well-defined
since $\alpha$ is uniquely determined.

On the remaining elements in the crystal, we specify the values of $\Psi$ as follows:
\begin{eqnarray}\label{negvalue}
&&\Psi(x_{-\theta}) = y_i \ot x_{-\theta} \\
&&\Psi(x_{-\gamma}) = x_\beta \ot x_{-\theta + \alpha} \qquad  \hbox{\rm if} \ \ \Psi(x_\gamma) = x_{\theta-\alpha} \ot x_{-\beta} \ \hbox{\rm and} \ \alpha \geq 0,  \nonumber  \\
&&\Psi(y_i) = x_{\theta - \alpha_i} \ot x_{-\theta + \alpha_i} \quad \  \hbox{\rm if vertex $i$ is connected to vertex $0$}, \nonumber \\
&&\Psi(y_j) =
x_{\theta-\alpha_{i_1} - \cdots - \alpha_{i_t}} \ot x_{-\theta + \alpha_{i_1} + \cdots + \alpha_{i_t}} \quad 
\hbox{\rm if} \ \ j \neq i,  \nonumber \end{eqnarray}
where  $\{i_1 = i, i_2, \dots, i_t = j\}$
is the unique (ordered) sequence of vertices  
connecting vertex $i$ and vertex $j$ in the Dynkin diagram.   
 
By its construction, the map $\Psi$ preserves the weight.
We claim that $\Psi$ commutes with the Kashiwara operators and so
is a crystal morphism.   We will check our assertion for $\tilde f_k$ ($k=1,\dots,n$) by considering
the various possibilities for $\gamma \in \Lambda_j^+$  \ $j=0,1,2$: \m

\noindent (Case $j=2$) \ \  If $\gamma = \theta$, then $\tilde f_i x_\theta = x_{\theta - \alpha_i}$
while $\tilde f_k x_\theta = 0$ for all $k \neq i$.   Then
$\tilde f_i \Psi(x_\theta) = \tilde f_i(x_\theta \ot y_i) = x_\theta \ot x_{-\alpha_i}$.
But $\Psi(\tilde f_i x_\theta) =  \Psi( x_{\theta-\alpha_i}) = x_\theta \ot x_{-\alpha_i}$ as $\theta-\alpha_i \in \Lambda_1^+$.   For $k \neq i$ ,   $\tilde f_k \Psi(x_\theta)  = 0 = \Psi(\tilde f_k x_\theta)$, so
$\Psi$ commutes with $\tilde f_k$ for all $k$ when applied to $x_\theta$.    \m

\noindent (Case $j=1$) \ \   We may assume  $\gamma = \theta - \beta \in \Lambda_1^+$ and
$\Psi(x_\gamma) = x_\theta \ot x_{-\beta}$.   If $\tilde f_k x_\gamma \neq 0$,
then $\gamma-\alpha_k \in \Lambda$, and if $k \neq i$, then $\gamma -\alpha_k
\in \Lambda_1^+$.    In that case, $\Psi(\tilde f_k x_\gamma) =
\Psi(x_{\gamma-\alpha_k}) = x_\theta \ot x_{-\beta-\alpha_k}$.    Since
$\tilde f_k \Psi(x_\gamma) = \tilde f_k (x_\theta \ot x_{-\beta}) = x_\theta \ot x_{-\beta -\alpha_k}$,
the result holds in this case.   Now when $k = i$ and $\tilde f_i x_\gamma \neq 0$, then $\gamma -\alpha_i \in \Lambda^+_0$.     The ``$\alpha$'' corresponding to $\gamma-\alpha_i$ is
$\alpha_i$ in this situation, and so  $\Psi(\tilde f_i x_\gamma) = \Psi(x_{\gamma-\alpha_i}) = x_{\theta-\alpha_i} \ot x_{-\beta}$.
Note that  $\tilde f_i \Psi(x_\gamma) = \tilde f_i(x_\theta \ot x_{-\beta}) = x_{\theta - \alpha_i} \ot x_{-\beta}$ because $-\beta + \alpha_i \not \in \Lambda$.   Indeed, if
$-\beta + \alpha_i \in \Lambda$, then $\theta = \gamma + \beta = (\gamma + \alpha_i)
+ (\beta -\alpha_i)$ implies that $\gamma + \alpha_i \in \Lambda^+$.    Moreover, since
$\gamma + \alpha_i \in \Lambda_2^+$, it must be that $\gamma + \alpha_i = \theta$.  But
then $\varphi_i(x_\theta) \geq 2$, a contradiction.
Consequently, $\Psi$ commutes with $\tilde f_i$ also, and the $j=1$ case is handled.
 \m

\noindent (Case $j=0$) \ \    We assume now that $\gamma \in \Lambda_0^+$
and let $\un j = (j_1,\dots, j_t = i)$ be the corresponding sequence.   Then
$\Psi(x_\gamma) = x_{\theta -\alpha} \ot x_{-\beta}$ (where
$\alpha = \alpha_{j_1} + \cdots + \alpha_{j_t}$),  \ $\gamma + \alpha \in
\Lambda_1^+$, and $\theta = \gamma + \alpha + \beta$.    If
$\tilde f_k x_\gamma \neq 0$, then $k \neq i$,  \ $k \in \hbox{\rm supp}(\gamma)$,
and $\gamma -\alpha_k \in \Lambda_0^+$.
  We consider two possible scenarios:

(1) Node  $k$ is connected to node $j_1$ in the Dynkin diagram of $\widehat \g$:
\   In this case,  $(\gamma -\alpha_k) + (\alpha + \alpha_k) \in \Lambda_1^+$,
and therefore $\Psi(\tilde f_k x_\gamma) = \Psi(x_{\gamma -\alpha_k})
= x_{\theta - (\alpha+\alpha_k)} \ot x_{-\beta}$.
On the other hand,   we claim that $\tilde f_k \Psi(x_\gamma) =
\tilde f_k (x_{\theta -\alpha} \ot x_{-\beta}) = x_{\theta-\alpha-\alpha_k} \ot x_{-\beta}$.
Indeed, $\theta = \gamma + \alpha + \beta =
(\gamma - \alpha_k) + (\alpha + (l+1)\alpha_k) + (\beta - l \alpha_k)$
whenever $\beta - l \alpha_k \in \Lambda^+$.   That is to say $\theta -\alpha - (l+1)\alpha_j \in
\Lambda$, which implies $\varphi_k(\theta-\alpha) \geq \varepsilon_k (-\beta) + 1$.
Thus $\tilde f_k$ should act on the first component of  
$x_{\theta -\alpha} \ot x_{-\beta}$ according to the tensor product rules;  hence
$\tilde f_k (x_{\theta -\alpha} \ot x_{-\beta}) = x_{\theta-\alpha-\alpha_k} \ot x_{-\beta}$
as asserted, and $\tilde f_k$ commutes with $\Psi$ when applied to $x_\gamma$.

(2) Node  $k$ is not connected to node $j_1$:  \   Then $\gamma - \alpha_k +\alpha \in
\Lambda_1^+$.  It follows from the $\Lambda_1^+$ case above that
$\Psi(\tilde f_k x_{\gamma}) = x_{\theta-\alpha} \ot x_{-\beta-\alpha_k}$.
Now $\varphi_k(\theta-\alpha) = 0$ since $k$ is not connected to $j_1$.    So
$\tilde f_k$ should act on the second component of $x_{\theta-\alpha} \ot x_{-\beta}$,
which implies $\tilde f_k\Psi(x_{\gamma}) =
\tilde f_k\left(x_{\theta-\alpha} \ot x_{-\beta}\right) = x_{\theta-\alpha} \ot x_{-\beta-\alpha_k}$.
Thus, $\tilde f_k$ commutes with $\Psi$ for such $x_\gamma$.     As this
is the last case that needed
to be considered, we have shown that $\tilde f_k$ commutes with $\Psi$ on all of $\mathcal B(\theta)$.

Now $\varphi_k$ and $\varepsilon_k$ have the same values on
$x_\theta$ and $x_\theta \ot y_i$ for each $k$.    Since $\Psi$
commutes with the Kashiwara operators, these functions will have
the same values on elements in the crystals $\mathcal B(\theta)$
and $\mathcal C(x_\theta \ot y_i)$ which correspond under $\Psi$,
because they are connected to $x_\theta$ and $x_\theta \ot y_i$
respectively by the same sequences of Kashiwara operators.  Hence $\Psi:\mathcal B(\theta) \iso \mathcal C(x_\theta \ot y_i)$ is a crystal isomorphism for all
types of affine algebras  except for $\hbox{\rm A}_n^{(1)}$ and
$\hbox{\rm C}_n^{(1)}$, 
which we now address. 

\m

\textbf {Case $\boldsymbol {\hbox{\rm A}_n^{(1)}}$:}     When $n = 1$, we have
$\theta = \alpha_1$, and 
\begin{eqnarray*}  x_{\alpha_1} &\mapsto& x_{\alpha_1} \ot y_1 \\
y_1 & \mapsto & y_1 \ot y_1 \\
x_{-\alpha_1} & \mapsto & y_1 \ot x_{-\alpha_1} \end{eqnarray*}
defines a crystal isomorphism in this case.  
 
For $n \geq 2$,   there are two vertices connected to the 0 vertex in the Dynkin
diagram, namely  the first  and  last.        When $i = 1$,  define a weight-preserving map
$\Psi: \mathcal B(\theta) \longrightarrow \mathcal C(x_\theta \ot y_1)$ by
assigning

\begin{eqnarray} x_{ \alpha_j+ \cdots + \alpha_k} &\mapsto& x_{\alpha_1 + \cdots + \alpha_k}
\ot x_{-\alpha_1 - \cdots-\alpha_{j-1}} \qquad \ (1 \leq j < k \leq n) \label{eq:Acase} \\
x_{\alpha_j} & \mapsto & x_{\alpha_1 + \cdots + \alpha_j}
\ot x_{-\alpha_1 - \cdots-\alpha_{j-1}} \qquad \ \, (2 \leq j \leq n) \nonumber \\
x_{\alpha_1} & \mapsto & x_{\alpha_1} \ot y_1 \nonumber \\
y_1 & \mapsto & y_1 \ot y_1 \nonumber \\
y_j& \mapsto  &
x_{\alpha_1 + \cdots + \alpha_{j-1}}
\ot x_{-\alpha_1 - \cdots-\alpha_{j-1}} \qquad (2 \leq j \leq n), \nonumber
\end{eqnarray}

\noindent and by using (\ref{negvalue})  for the negative roots. 
One can check that the Kashiwara operators $\tilde f_l$ commute with $\Psi$
and proceed as before to show that $\Psi$ is a crystal isomorphism.  For example,
in the case of $x_{ \alpha_j + \cdots + \alpha_k}$,
there are exactly two Kashiwara operators $\tilde f_j$ and $\tilde f_k$ whose action
on each side of \eqref{eq:Acase}
is nonzero.        For the first,  
\begin{eqnarray*} && \Psi(\tilde f_j x_{ \alpha_j + \cdots + \alpha_k})
= x_{\alpha_{j+1} + \cdots + \alpha_k} \mapsto  x_{\alpha_1 + \dots + \alpha_k} \ot
x_{-\alpha_1 - \cdots - \alpha_j},  \quad \hbox{\rm while}  \\
&& \tilde f_j\Psi(x_{ \alpha_j + \cdots + \alpha_k}) =
\tilde f_j \left (x_{\alpha_1 + \cdots + \alpha_k}
\ot x_{-\alpha_1 - \cdots-\alpha_{j-1}} \right) =
x_{\alpha_1 + \cdots + \alpha_k} \ot x_{-\alpha_1 - \cdots - \alpha_j}.  \end{eqnarray*}
For the second,
\begin{eqnarray*}&&\Psi(\tilde f_kx_{ \alpha_j + \cdots + \alpha_k})
= x_{\alpha_j + \cdots + \alpha_{k-1}} \mapsto  x_{\alpha_1 + \dots + \alpha_{k-1}} \ot
x_{-\alpha_1 - \cdots - \alpha_{j-1}}, \quad \hbox{\rm while} \\
&&\tilde f_k\Psi(x_{ \alpha_j+ \cdots + \alpha_k}) =
\tilde f_k\left (x_{\alpha_1 + \cdots + \alpha_k}
\ot x_{-\alpha_1 - \cdots-\alpha_{j-1}} \right) =
x_{\alpha_1 + \cdots + \alpha_{k-1}} \ot x_{-\alpha_1 - \cdots - \alpha_{j-1}}. 
\end{eqnarray*} 
The remaining calculations to verify that $\Psi$ commutes with the Kashiwara
operators and is a crystal isomorphism are similar and are left to the reader.

Assume now that  $i = n$  and define a weight-preserving map
$\Psi: \mathcal B(\theta) \longrightarrow \mathcal C(x_\theta \ot y_n)$ by
assigning

\begin{eqnarray} x_{ \alpha_j + \cdots + \alpha_k} &\mapsto& x_{\alpha_j + \cdots + \alpha_n}
\ot x_{-\alpha_{k+1} - \cdots-\alpha_{n}} \qquad \ (1 \leq j < k \leq n) \label{eq:Acasen} \\
x_{\alpha_j} & \mapsto & x_{\alpha_j+ \cdots + \alpha_n}
\ot x_{-\alpha_{j+1} - \cdots-\alpha_{n}} \qquad \ \, (2 \leq j \leq n) \nonumber \\
x_{\alpha_n} & \mapsto & x_{\alpha_n} \ot y_n\nonumber \\
y_j & \mapsto  &
x_{\alpha_{j+1} + \cdots + \alpha_{n}}
\ot x_{-\alpha_{j+1} - \cdots-\alpha_{n}} \qquad (2 \leq j \leq n)  \nonumber \\
y_n & \mapsto & y_n \ot y_n. \nonumber
\end{eqnarray}    The argument is exactly as before.  For example,  the computations involving the first expression in \eqref{eq:Acasen}
are these:

\begin{eqnarray*} \Psi(\tilde f_j x_{ \alpha_j + \cdots + \alpha_k} ) &=&
\Psi(x_{ \alpha_{j+1}+ \cdots + \alpha_k}) =  x_{\alpha_{j+1} + \cdots + \alpha_n}
\ot x_{-\alpha_{k+1} - \cdots-\alpha_{n}}\\
\tilde f_j \Psi(x_{ \alpha_j + \cdots + \alpha_k}) &=&
\tilde f_j \left(x_{\alpha_j + \cdots + \alpha_n}
\ot x_{-\alpha_{k+1} - \cdots-\alpha_{n}}\right) =
x_{\alpha_{j+1} + \cdots + \alpha_n}
\ot x_{-\alpha_{k+1} - \cdots-\alpha_{n}} \\
\Psi(\tilde f_k x_{ \alpha_j+ \cdots + \alpha_k} ) &=&
\Psi(x_{ \alpha_j+ \cdots + \alpha_{k-1}}) =  x_{\alpha_j + \cdots + \alpha_n}
\ot x_{-\alpha_k - \cdots-\alpha_{n}}\\
\tilde f_k\Psi(x_{ \alpha_j + \cdots + \alpha_k}) &=&
\tilde f_k \left(x_{\alpha_j + \cdots + \alpha_n}
\ot x_{-\alpha_{k+1} - \cdots-\alpha_{n}}\right) =
x_{\alpha_j + \cdots + \alpha_n}
\ot x_{-\alpha_k - \cdots-\alpha_{n}}. \end{eqnarray*}
\m

\textbf {Case $\boldsymbol {\hbox{\rm C}_n^{(1)}}$:}   In this case, $\theta =
2 \alpha_1 + 2 \alpha_2 + \cdots + 2 \alpha_{n-1} + \alpha _n$ and vertex 0 is
connect to vertex 1.   We define a weight-preserving map
$\Psi: \mathcal B(\theta) \longrightarrow \mathcal C(x_\theta \ot y_1)$ by
assigning

\begin{eqnarray} &&\hspace {3truein} \label{eq:Ccasen} \\
\vspace{-.2 truein}
x_{ \alpha_j + \cdots + \alpha_{k-1} + 2\alpha_k + \cdots + 2\alpha_{n-1} + \alpha_n} &\mapsto& x_{ \alpha_1 + \cdots + \alpha_{k-1} + 2\alpha_k + \cdots + 2\alpha_{n-1} + \alpha_n}
\ot x_{-\alpha_{1} - \cdots-\alpha_{j-1}} \nonumber \\
&& \hspace{1.8truein}  (1 \leq j \leq k \leq n-1) \nonumber  \\
x_{\alpha_j + \cdots + \alpha_k} & \mapsto
& x_{\alpha_1+ \cdots + \alpha_k} \ot x_{-\alpha_1 - \cdots - \alpha_{j-1}} \qquad (1 \leq j \leq k \leq n) \nonumber \\ x_{\alpha_1} & \mapsto & x_{\alpha_1} \ot y_1\nonumber \\
x_{\alpha_j} & \mapsto & x_{\alpha_1 +\cdots + \alpha_{j}} \ot  x_{-\alpha_1 - \cdots - \alpha_{j-1}}
\qquad \quad \ (2 \leq j \leq n) \nonumber \\
y_1 & \mapsto & y_1 \ot y_1. \nonumber \\
y_j & \mapsto  &
x_{\alpha_1 + \cdots + \alpha_j}
\ot x_{-\alpha_1 - \cdots-\alpha_{j}}\qquad \qquad (2 \leq j \leq n)  \nonumber
\end{eqnarray}
and using  (\ref{negvalue}).   The arguments are  straightforward as in the previous cases.  Here is a sample computation
for $x_{\gamma}$,  where $\gamma = \alpha_j + \cdots + \alpha_{k-1} + 2\alpha_k + \cdots + 2\alpha_{n-1} + \alpha_n$:

\begin{eqnarray*} \Psi(\tilde f_j x_{\gamma}) &=&\Psi(x_{\gamma-\alpha_j} ) \\
& =&x_{\alpha_1 + \cdots + \alpha_{k-1} +2 \alpha_k + \cdots + 2\alpha_{n-1} + \alpha_n}
\ot x_{-\alpha_{1} - \cdots-\alpha_j}  \\
\tilde f_j \Psi(x_{\gamma}) &=&
\tilde f_j \left(x_{\alpha_1 + \cdots + \alpha_{k-1} + 2\alpha_k + \cdots + 2\alpha_{n-1} + \alpha_n}
\ot x_{-\alpha_{1} - \cdots-\alpha_{j-1}}\right) \\
&=&
x_{\alpha_1 + \cdots + \alpha_{k-1} + 2\alpha_k + \cdots + 2\alpha_{n-1} + \alpha_n}
\ot x_{-\alpha_{1} - \cdots-\alpha_{j}}  \\
\Psi(\tilde f_k x_{ \gamma}) &=&
\Psi(x_{ \gamma-\alpha_k}  ) \\ &=&  x_{\alpha_1 + \cdots + \alpha_{k-1} + \alpha_k + 2\alpha_{k+1}+ \cdots + 2\alpha_{n-1} + \alpha_n} \ot x_{-\alpha_{1} - \cdots-\alpha_{j-1}} \\
\tilde f_k\Psi(x_{ \gamma}) &=&
\tilde f_k \left(x_{\alpha_1 + \cdots + \alpha_{k-1} +2 \alpha_k + \cdots + 2\alpha_{n-1} + \alpha_n}
\ot x_{-\alpha_{1} - \cdots-\alpha_{j-1}}\right) \\ &=&
x_{\alpha_1 + \cdots + \alpha_{k-1} + \alpha_k+ 2 \alpha_{k+1} + \cdots + 2\alpha_{n-1} + \alpha_n}
\ot x_{-\alpha_{1} - \cdots-\alpha_{j-1}}.   \end{eqnarray*}
\end{proof}
\m

\begin{ex}   When $\widehat \g = \hbox{\rm D}_4^{(3)}$ and $\g = \hbox{\rm G}_2$,
then $\mathcal B(\theta)$ is the crystal of the 7-dimensional irreducible $U_q(\g)$-module,
and the crystal isomorphism  $\mathcal C(x_\theta \otimes y_1) \cong \mathcal B(\theta)$  
gives   the  ``crystal  octonion algebra''  (without the unit element).
Displayed below is a portion of  the multiplication table in this algebra.
All other products  between basis elements are zero.   

\begin{equation*}  \vbox{
\offinterlineskip \halign{ \strut \vrule  $\hfil#\hfil$  & \vrule \vrule
$\hfil#\hfil$  & \vrule $\hfil#\hfil$ & \vrule $\hfil#\hfil$ & \vrule $\hfil#\hfil$  
\vrule\cr \noalign{\hrule} \, 
& x_{-\theta} & x_{-\alpha_1-\alpha_2} & x_{-\alpha_1}
 & y_1   \cr \noalign{\hrule}  \noalign{\hrule} \,
x_\theta &\, 0 & \,x_{\alpha_1} & \, x_{\alpha_1+\alpha_2} & \, x_\theta   \cr \noalign{\hrule}
 \,
x_{\alpha_1+\alpha_2} & \, x_{-\alpha_1} & \, y_1 & \, 0 & \, 0 \cr \noalign {\hrule}\, 
x_{\alpha_1} & \, x_{-\alpha_1-\alpha_2} & \, 0 & \, 0 & \, 0  \cr \noalign{\hrule}  \,
y_1 & \, x_{-\theta} &\, 0 & \, 0 & \, 0   \cr \noalign{\hrule}} }  \end{equation*}

\end{ex} 
\m

\begin{section} {The Energy Function}  \end{section}

The energy function on the perfect crystal  $\mathcal B$  enables us to determine the weight of
a path  (see Theorem \ref{thm:wtchar}).    Here we compute the energy function, which
we assume has been normalized so that $H(\emptyset \ot \emptyset) = 0$.   First, observe
how the maximal vectors of the right side of  \eqref{eq:dec}  are connected:

\begin{equation}\label{eq:energytable} \begin{array}{ccccc}
x_\theta \ot x_\theta &\zlarrow& x_\theta \ot \emptyset &\zlarrow& x_\theta \ot x_{-\theta}   \\
&&&&  \uparrow_0 \\
&&&& \emptyset \ot x_{-\theta} \\
x_\theta \ot x_{\theta-\alpha} &\zlarrow& x_\theta \ot x_{-\alpha} &&\ \  \uparrow_{\ast \dots \ast} \\
&& \uparrow_0 &&    \\
&& \emptyset \ot x_{-\alpha} &\dotlarrow& \emptyset \ot  x_\theta \\
&&& {}_{\ast \dots \ast}\swarrow &\, \uparrow_0 \\
x_\theta \ot y_i&\zlarrow& \emptyset \ot y_i && \emptyset \ot \emptyset .
\end{array}\end{equation} \m

Thus, by \eqref{eq:ef} we deduce

\begin{equation}\label{eq:efvalues} \vbox{
\offinterlineskip \halign{ \strut \vrule  $\hfil#\hfil$  & \vrule $\hfil#\hfil$  & \vrule
$\hfil#\hfil$  & \vrule $\hfil#\hfil$ & \vrule $\hfil#\hfil$ & \vrule $\hfil#\hfil$ & \vrule  $\hfil#\hfil$  & \vrule  $\hfil#\hfil$
\vrule\cr \noalign{\hrule} \,\hbox{\rm max'l vec.} &\, \emptyset \ot \emptyset & \, \emptyset \ot x_\theta & \, x_\theta \ot x_{-\theta}  & \,
x_{\theta} \ot \emptyset
 & \,x_\theta \ot x_\theta  & \,x_\theta \ot x_{\theta -\alpha}  &\, x_\theta \ot y_i  \cr \noalign{\hrule}\,
\hbox{\rm $H$-value} &\, 0 &\, 1 & \, 0 & \, 1 & \,  2 & \, 1
& \, 0  \cr \noalign{\hrule}} }\end{equation}    \bi

Since the values of $H$ are constant on each connected component of  the crystal graph $\mathcal B \ot \mathcal B$
without the 0-arrows,  we need to describe the connected components containing the maximal
vectors above.
\m

\begin{equation}  \label{eq:cc} \end{equation}
\begin{itemize}
\item[{\rm  (1)}]   The connected component \ $\mathcal C(x_\theta \ot x_{-\theta})$ \ containing the
maximal vector \ $x_\theta \ot x_{-\theta}$ is a singleton $\{x_\theta \ot x_{-\theta}\}$ and by
the above table,  $H = 0$.  \m

\item[{\rm  (2)}]  Similarly,  $\mathcal C(\emptyset \ot \emptyset) = \{\emptyset \ot \emptyset\}$
and $H = 0$.  \m

\item[{\rm  (3)}]  The connected components containing $x_\theta \ot \emptyset$ and
$\emptyset \ot x_\theta$ have a simple description,

\begin{eqnarray*} \mathcal C(x_\theta \ot \emptyset) &=& \{ b \ot \emptyset \mid
b \in \mathcal B(\theta) \}\\
\mathcal C(\emptyset \ot x_\theta) &=& \{ \emptyset \ot b \mid b \in \mathcal B(\theta)\} \end{eqnarray*}
and on both, $H = 1$.  \m
\end{itemize}

Next,   we will describe the component $\mathcal C(x_\theta \ot x_{\theta})$.
The components of the form  $\mathcal C(x_\theta \ot x_{\theta -\alpha})$ will then
consist of all the tensor products of elements of $\mathcal B$
which have not been already specified,  and on them, $H = 1$.

\bi
\begin{prop}\label{2thetacomp} If $\widehat \g \neq \hbox{\rm A}_{2n}^{(2)}$, the connected component of the maximal vector $x_\theta \ot x_\theta$
is given by

\begin{eqnarray*}
\mathcal C(x_\theta \ot x_\theta) &=& \{ x_\alpha \ot x_\beta \mid \alpha,\beta \in \Lambda^+\sqcup \Lambda^-, \ \alpha \leq \beta\} \\
&&\textstyle{\bigsqcup} \, \{ y_i \ot x_\beta \mid \theta(h_i) > 0, \ \beta=\alpha_i+\gamma \in \Lambda^+
\ \hbox{\rm for some}\ \gamma \in \Lambda\} \\
&& \textstyle{ \bigsqcup} \, \{ x_{-\alpha} \ot y_i \mid \theta(h_i) > 0, \ \alpha=\alpha_i+\gamma \in \Lambda^+\ \hbox{\rm for some}\ \gamma \in \Lambda\},
\end{eqnarray*}
where $\alpha \leq \beta$ if and only if $\beta-\alpha \in
\sum_{i\in \mathcal I\setminus \{0\}} \mathbb Z_{\geq 0}  \alpha_i$. \m

If $\widehat \g = \hbox{\rm A}_{2n}^{(2)}$, we have
\begin{equation*}
\mathcal C(x_\theta \ot x_\theta) = \{ x_\alpha \ot x_\beta \mid
\alpha,\beta \in \Lambda^+\sqcup \Lambda^-, \ \alpha \leq \beta\}.
\end{equation*}

\end{prop}

\begin{proof}

\vspace{.3 truein}

\noindent \textbf{Step 1.}  \, \ {\it For any  $\alpha \in \Lambda^+ \sqcup \Lambda^-$,  there exist   sequences $\un i = (i_1,\dots,i_r)$ and $\un j = (j_1,\dots, j_t)$ of indices in $\mathcal I \setminus \{0\}$  such that
$x_\theta  \ot x_\theta \iunarrow x_\alpha \ot x_\theta$ and
$x_\theta  \ot x_\theta {\buildrel {\un  j} \over \rightarrow} x_\alpha \ot x_\alpha$.   } \m

 (Here the shorthand $\iunarrow $ signifies  that  the sequence of
arrows
${\buildrel {i_1} \over \rightarrow} \  {\buildrel {i_2} \over \rightarrow}  \ \cdots \ {\buildrel {i_r} \over \rightarrow}$ has been applied.)

 By the tensor product rule,  the Kashiwara
 operator $\tilde e_{k}$ acts on the second component
 first, so using that and  the connectedness of  $\mathcal B(\theta)$,
 we see that both $x_{\alpha} \ot x_{\theta}$
 and $x_{\alpha} \ot x_{\alpha}$ can be connected to $x_\theta \ot x_{\theta}$
 by applying a suitable sequence of operators $\tilde e_{k}$.

 \m

\noindent \textbf{Step 2.} \,  {\it Suppose that $\alpha, \beta \in \Lambda^+\sqcup \Lambda^-$ and $\alpha \leq \beta$.   Then there exist  $\gamma \in \Lambda$ and a sequence $\un j$  such that}  $x_\gamma \ot x_\theta\,
{\buildrel {\un j} \over \rightarrow}\, x_\alpha \ot x_\beta$.    Hence $x_\alpha \ot x_\beta
\in \mathcal C(x_\theta \ot x_\theta)$.

\m

Working in  the crystal base  $\mathcal B(\theta)$,  we may choose sequences  $\un i_\alpha$ and $\un i_\beta$ such that
$x_\theta \, {\buildrel {\un i_\alpha} \over \rightarrow} \, x_\alpha$ and
$x_\theta \, {\buildrel{\un i_\beta} \over \rightarrow} \, x_\beta$.  We induct on  $r = |\, \un i_\alpha\,|$
and $s = |\,\un i_\beta\,|$.
To begin, if $r = 0$, then $\theta = \alpha \leq \beta \leq \theta$, which  forces all of them to be equal.
Thus, $x_\alpha \ot x_\beta = x_\theta \ot x_\theta$ in this case, and we may take
$\gamma = \theta$ and $\un j$  to be the empty sequence.
Similarly, if $s = 0$, then $\beta = \theta$ and $x_\alpha \ot x_\beta = x_\alpha \ot x_\theta$.
so the result follows from Step 1.

Suppose then  that $r > 0$ and $s > 0$.    Now if there exists an $\ell  \in \un i_\beta$ such that
$\beta + \alpha_\ell \in \Lambda$, but $\alpha -\alpha_\ell \not \in \Lambda$, then
$\tilde e_\ell (x_\alpha \ot x_\beta) = x_\alpha \ot x_{\beta + \alpha_\ell}$.
 By induction,
there exist $\gamma \in \Lambda$ and  a sequence $\un {j}'$ such that   $x_\gamma \ot x_\theta \,
{\buildrel {\un j'}  \over \rightarrow} \, x_\alpha \ot x_{\beta + \alpha_\ell}$.  Adjoining $\ell$ to the end
of the sequence ${\un j'}$   gives the desired sequence $\un j$.
Since by Step 1, there exists a sequence $\un i$  such that $x_\theta \ot x_\theta \iunarrow x_\gamma \ot x_\theta \, {\buildrel {\un j} \over \rightarrow}\, x_\alpha \ot x_\beta$, it follows  that $x_\alpha \ot x_\beta \in
 \mathcal C(x_\theta \ot x_\theta)$.

Suppose there exists a value of $k$ such that $\alpha + \alpha_k \in \Lambda$ and
$\alpha+\alpha_k \leq \beta$.  Then
$\tilde e_k(x_\alpha \ot x_\beta)  = x_{\alpha+\alpha_k} \ot x_\beta$.   By induction, there exist $x_\gamma$ with
$\gamma \in \Lambda$ and a sequence ${\un j'}$, such that
$x_\gamma \ot x_\theta {\buildrel {\un j'}  \over \rightarrow} \, x_{\alpha+\alpha_k} \ot x_{\beta}$.
Adjoining $k$ to the end of the sequence and using Step 1 then gives  the result.

Hence we may suppose that no such $\ell$ and $k$ exist.  Therefore,
(1) for all $k \in \un i_\alpha$, we have $\alpha+\alpha_k \not \in \Lambda$, i.e. $\tilde e_k x_\alpha = 0$;
and (2) for all $\ell \in \un i_\beta$ with $\beta +\alpha_\ell \in \Lambda$, we have
$\alpha-\alpha_\ell \not \in \Lambda$.   By the tensor product rule, this implies that
$\tilde e_m(x_\alpha \ot x_\beta) = 0$ for all $m \in \mathcal I \setminus \{0\}$.
Hence $x_\alpha \ot x_\beta$ is a maximal vector in $\mathcal B(\theta) \ot \mathcal B(\theta)$.
But by \cite[Cor.~4.4.4\,(2)]{HK}, it must be that $x_\alpha$ is a maximal vector of $\mathcal B(\theta)$.
Thus, $\theta = \alpha \leq \beta \leq \theta$, which forces  $r=0=s$.  Since we are
assuming both $r$ and $s$ are positive, this case cannot happen.   The proof of Step 2 is now
complete. \m

\noindent \textbf{Step 3.} \, {\it The elements $y_i \ot x_\beta$ with $\theta(h_i) > 0$
and $\beta = \alpha_i + \gamma \in \Lambda^+$ for some $\gamma \in \Lambda$ belong to}
$\mathcal C(x_\theta \ot x_\theta)$. \m

The assumptions imply that  $\tilde e_i(y_i \ot x_\beta) = x_{\alpha_i} \ot x_\beta$.   Then by
Steps 1 and 2, we can find $\gamma \in \Lambda$ and sequences $\un i$ and $\un j$
so that
$x_\theta \ot x_\theta \iunarrow x_\gamma \ot x_\theta\,
{\buildrel {\un j} \over \rightarrow}\, x_{\alpha_i} \ot x_\beta \, \iarrow \, y_i \ot x_\beta.$
Hence, $y_i \ot x_\beta \in \mathcal C(x_\theta \ot x_\theta)$, as claimed.

\m

\noindent \textbf{Step 4.} \,  {\it  The elements $x_{-\alpha} \ot y_i$ such that
$\theta(h_i) > 0$ and  $\alpha = \alpha_i + \gamma
\in \Lambda^+$ for some $\gamma \in \Lambda$ belong to  $ \mathcal C(x_\theta \ot x_\theta)$.}
\m

Observe here that because $-\alpha -\alpha_i \in \Lambda^-$,
$\tilde f_i (x_{-\alpha} \ot y_i) = x_{-\alpha} \ot  x_{-\alpha_i}$,  where $-\alpha \leq  -\alpha_i$.
Thus, by Step 2, we can find a sequence $\un j$ such that
$x_\gamma \ot x_\theta \,{\buildrel {\un j} \over \rightarrow}\,  x_{-\alpha}  \ot x_{-\alpha_i} \iarrow
x_{-\alpha} \ot y_i$, which implies that $x_{-\alpha} \ot y_i \in  \mathcal C(x_\theta \ot x_\theta)$
by Step 1. \m

We have argued in Steps 1-4 that the sets listed in the statement
of the lemma belong to the connected component $\mathcal
C(x_\theta \ot x_\theta)$. It remains to show that the set in the
right-hand side of our statement is closed under the action of
Kashiwara operators, which can be done in a rather straightforward
way using the tensor product rule. For instance, suppose $\alpha,
\beta \in \Lambda$ and $\alpha \le \beta$. Then by the tensor
product rule, we have $\tilde f_i (x_{\alpha} \ot x_{\beta}) =
x_{\alpha - \alpha_i} \ot x_{\beta}$ or $\tilde f_i (x_{\alpha}
\ot x_{\beta}) = x_{\alpha} \ot x_{\beta - \alpha_i}$. In the
first case, we are done. In the latter case, since $\alpha \le
\beta$, we have $\un i_\beta \subset \un i_\alpha$ and $i \in \un
i_\alpha \setminus \un i_\beta$. Hence $\alpha \le \beta -
\alpha_i$.  The other cases can be verified similarly.
 
 \end{proof}
 
 \medskip
 Combining the results of the last two sections, we arrive at the following:
 \bigskip
 \begin{thm}\label{thm:energy}   The energy function $H$ on the level 1 perfect crystal $\mathcal B$ is given explicitly by Table \eqref{eq:efvalues},   Proposition \ref{thetacomp}, and Proposition \ref{2thetacomp}.
 \end{thm}
 
 \medskip
 The energy function formula  for the weight of a path in Theorem \ref{thm:wtchar}
 (in particular,  the coefficient of $\delta$)  gives  the homogeneous
 grading in the path realization of the crystal base of any level 1 (basic)  representation of
 $U_q(\widehat \g)$.      Theorem  \ref{thm:energy} and the general character formula in
 Theorem \ref{thm:wtchar} then yield new character formulas for the basic representations.  Their
 explicit combinatorial expressions in terms of the roots deserve further study.

\end{document}